\documentclass[12pt, a4paper]{article}
\usepackage{amsmath, amssymb, amsthm, amscd}

\newtheorem{theorem}{Theorem}[section]

\newtheorem{lemma}[theorem]{Lemma}
\newtheorem{proposition}[theorem]{Proposition}
\newtheorem {corollary}[theorem]{Corollary}
\theoremstyle {definition}
\newtheorem {definition}[theorem]{Definition}
\newtheorem {example}[theorem]{Example}
\theoremstyle {remark}
\newtheorem{remark}[theorem]{Remark}

\def\Spec{\operatorname{Spec}}
\def\Ann{\operatorname{Ann}}
\def\Ass{\operatorname{Ass}}
\def\Supp{\operatorname{Supp}}
\def\NCM{\operatorname{NCM}}

\def\Rad{\operatorname{Rad}}

\newcommand{\fm}{\ensuremath{\mathfrak m}}
\newcommand{\fa}{\ensuremath{\mathfrak a}}
\newcommand{\fp}{\ensuremath{\mathfrak p}}
\newcommand{\fq}{\ensuremath{\mathfrak q}}
\newcommand{\un}{\ensuremath{\underline}}
\newcommand{\F}{\ensuremath{\mathcal F}}
\newcommand{\D}{\ensuremath{\mathcal D}}

\newcommand{\ha}{\ensuremath{I_{\mathcal F,M}}}

\begin{document}
\title{On the Structure of Sequentially Generalized Cohen-Macaulay Modules}
\author{Nguyen Tu Cuong\footnote{Email: ntcuong@math.ac.vn} \ and Doan Trung Cuong\footnote{Email: dtcuong@math.ac.vn}\\
Institute of Mathematics,\\
18 Hoang Quoc Viet Road, 10307 Hanoi, Vietnam}
\date{ }
\maketitle
\begin{abstract}
A finitely generated module $M$ over a local ring is called a sequentially generalized Cohen-Macaulay module if there is a filtration of submodules of $M$:\ $M_0\subset M_1\subset \ldots \subset M_t=M$ such that $\dim M_0<\dim M_1< \ldots <\dim M_t$ and each $M_i/M_{i-1}$ is generalized Cohen-Macaulay. The aim of this paper is to study the structure of this class of modules. Many basic properties of these modules are presented and various characterizations of sequentially generalized Cohen-Macaulay property by using local cohomology modules, theory of multiplicity and in terms of systems of parameters are given. We also show that the notion of dd-sequences defined in \cite{cc} is an important tool for studying this class of modules. \vspace{.2cm}
\end{abstract}

\section{Introduction}
Let $(R, \fm)$ be a commutative Noetherian local ring and $M$ a finitely generated $R$-module of dimension $d$. Let $\un x=(x_1,\ldots, x_d)$ be a system of parameters of $M$. It is well-known that the length $\ell(M/\un xM)$ carries a lot of information about the structure of $M$. If $\ell(M/\un xM)=e(\un x;M)$, where $e(\un x;M)$ is the Serre multiplicity of $M$ relative to $\un x$, then $M$ is a Cohen-Macaulay module. The notion of Buchsbaum modules introduced by St\"uckrad and Vogel is the first extension of Cohen-Macaulay modules, it contains all modules such that the difference $\ell(M/\un xM)-e(\un x;M)$ is a constant for all systems of parameters $\un x$. A further generalization was obtained by Schenzel, Trung and the first author in \cite{cst}, they considered the class of modules $M$ such that for all systems of parameters $\un x$ the difference $\ell(M/\un xM)-e(\un x;M)$ is bounded above by a constant. This is equivalent to the fact that there is a system of parameters $\un x$ such that $\ell(M/(x_1^{n_1}, \ldots, x_d^{n_d})M)=n_1\ldots n_de(\un x; M)+c$ for all $n_1, \ldots, n_d>0$, where $c$ is a constant. These modules have many similar properties as of Cohen-Macaulay modules and were called generalized Cohen-Macaulay modules. The theory of generalized Cohen-Macaulay modules was developed rapidly in the 1980's and early 1990's by the works of many authors and found its applications in many fields of commutative algebra and algebraic geometry. Another generalization of Cohen-Macaulay module is the notion of sequentially Cohen-Macaulay modules introduced first  by Stanley  \cite{st}. A module $M$ is called a sequentially Cohen-Macaulay module if there is a filtration $M_0 \subset M_1 \subset \ldots \subset M_t=M$ of submodules of $M$ such that each $M_{i+1}/M_i$ is Cohen-Macaulay and $\dim M_0<\dim M_1<\ldots <\dim M_t$. Historically, Stanley defined this notion for graded modules in order to study the so-called Stanley-Reisner rings (see also Herzog-Sbarra \cite{hs}). After that, this notion was defined for modules over local rings by Schenzel \cite{sch}, Nhan and the first author \cite{cn}. In the same paper, the  authors also introduced the notion of sequentially generalized Cohen-Macaulay module and gave a characterization for these modules in terms of local cohomology modules. The definition of sequentially generalized Cohen-Macaulay module is similar to the one of sequentially Cohen-Macaulay module except each module $M_{i+1}/M_i$ is required to be a generalized Cohen-Macaulay module instead of being Cohen-Macaulay. In this case, that a filtration is called a {\it generalized Cohen-Macaulay filtration}. The aim of this paper is to study basic properties of these modules with further purpose  toward a theory of sequentially generalized Cohen-Macaulay modules.

In order to study sequentially generalized Cohen-Macaulay modules, we consider a filtration $\F:\ M_0\subset M_1\subset \ldots \subset M_t=M$ of submodules of $M$, which satisfies the condition that  $\dim M_0<\dim M_1<\ldots <\dim M_t=d$. The most  important example of filtration satisfying the dimension condition as above is the dimension filtration. We say that a filtration $\F:\ M_0\subset M_1\subset \ldots \subset M_t=M$ is the {\it dimension filtration} of $M$ if each $M_i$ is the biggest submodule of $M_{i+1}$ with $\dim M_i<\dim M_{i+1}$ (cf. \cite{cn}, \cite{st}). For a filtration satisfying the dimension condition $\F$ with $d_i=\dim M_i$, we restrict ourself to those systems of parameters $\un x=(x_1, \ldots, x_d)$, which are called {\it good systems of parameters} of $M$, such that $M_i\cap (x_{d_i+1}, \ldots, x_d)M=0,\ i=0, 1, \ldots, t-1$. Then $(x_1,\ldots, x_{d_i})$ is a system of parameters of $M_i$. It is proved in \cite{cc1} that the difference 
$$I_{\F, M}(\un x)=\ell(M/\un xM)-\sum_{i=0}^te(x_1,\ldots, x_{d_i};M_i),$$
is a non-negative integer. From our point of view, $I_{\F, M}(\un x)$ is suitable to the study of sequentially Cohen-Macaulay and sequentially generalized Cohen-Macaulay modules. It has been shown by the authors in \cite{cc1} that $M$ is a sequentially Cohen-Macaulay module if and only if there is a filtration $\F$ and a good system of parameters  such that $$I_{\F, M}(x_1^{n_1}, \ldots, x_d^{n_d})=0, \text{ for all } n_1, \ldots, n_d>0,$$
or equivalently, $\ell(M/(x_1^{n_1}, \ldots, x_d^{n_d})M)=\sum_{i=0}^tn_1\ldots n_{d_i}e(x_1, \ldots, x_{d_i}; M_i)$. As one of the main results of this paper, we will show that $M$ is a sequentially generalized Cohen-Macaulay module if and only if there are a filtration $\F$ and a good system of parameters $\un x$ such that $I_{\F, M}(x_1^{n_1}, \ldots, x_d^{n_d})$ is a constant for all $n_1, \ldots, n_d>0$. Moreover, this constant is independent of the choice of systems of parameters and can be expressed in terms of length of certain local cohomology modules. The key in the proof of these results is the use of the notion of dd-sequence developed in \cite{cc}. dd-Sequence was first invented for a different purpose, see \cite{acta}, \cite{ta1}, \cite{ta2}. However, when studying the two classes of sequentially Cohen-Macaulay and sequentially generalized Cohen-Macaulay modules we found that this notion is very useful since all these modules admit such a sequence.

The paper is organized as follows.

In Section 2 we recall briefly some facts about filtrations satisfying the dimension condition and good systems of parameters. Some properties of dd-sequence defined in \cite{cc} are presented in this section.

In Section 3 we introduce the notion of generalized Cohen-Macaulay filtrations to investigate the structure  of sequentially generalized Cohen-Macaulay modules. We first show some properties of these modules by using local cohomology modules, localization, passing to quotient, etc. As the main result of this section, we show that for a sequentially generalized Cohen-Macaulay module $M$ there are a filtration $M_0\subset M_1\subset \ldots \subset M_t=M$ and a system of parameters $\un x=(x_1, \ldots, x_d)$ such that 
\begin{equation*}
\tag{*}\ell(M/(x_1^{n_1}, \ldots, x_d^{n_d})M)=\sum_{i=0}^tn_1\ldots n_{d_i}e(x_1, \ldots, x_{d_i}; M_i)+C
\end{equation*}
for all $n_1, \ldots, n_d>0$, where $d_i=\dim M_i$ and $C$ is a constant (Theorem \ref{p-stand}). 

We use Section 4 to study the constant $C$ in the equality (*). This number is  important in our investigation because it is the least bound for the function $\ha(x_1^{n_1}, \ldots, x_d^{n_d})$. The main result of this section is an expression of $C$ in terms of lengths of certain local cohomology modules, 
$$C=\sum_{i=0}^{t-1}\sum_{j=0}^{d_{i+1}-1}c_{ij}\ell(H^j_\fm(M/M_i),$$
where $c_{ij}=\sum_{k=d_i}^{d_{i+1}-1}\binom{k-1}{j-1}$. 

Using the theory of multiplicity we prove in Section 5 various characterizations of sequentially generalized Cohen-Macaulay modules in terms of good systems of parameters. Note that the filtrations of submodules of $M$ considered in this section are not necessary to be generalized Cohen-Macaulay filtrations.

The last section is devoted to study the Hilbert-Samuel function of a sequentially generalized Cohen-Macaulay module with respect to an ideal generated by a good system of parameters satisfying the equality (*). We compute all the coefficients of  the Hilbert-Samuel polynomial explicitely by using local cohomology modules.


\section{Preliminary}
Throughout this paper, $(R,\fm)$ is a commutative Noetherian local ring and $M$ is a finitely generated $R$-module of dimension $d$.

In this section we will recall briefly some basic facts about filtrations satisfying the dimension condition,  good systems of parameters defined in \cite{cc1}. Some preparations on dd-sequences and generalized Cohen-Macaulay modules are also presented.

\begin{definition}
(1) We say that a finite filtration of submodules of $M$
$$\mathcal F : \ M_0\subset M_1\subset \cdots \subset M_t=M$$
 satisfies the {\it  dimension condition} if
$\dim M_0<\dim M_1<\ldots <\dim M_{t-1}<\dim M$, and we also say in this case that the filtration $\mathcal F$ has the length $t$.
\item(2) A filtration $\D : D_0\subset D_1\subset \cdots \subset D_t=M$ is called the {\it dimension filtration} of $M$ if the following two conditions are satisfied

a) $D_{i-1}$ is the largest submodule of $D_i$ with $\dim D_{i-1}< \dim D_i$ for $i= t, t-1, \ldots, 1$;

b) $D_0=H_\fm^0(M)$ is the $0^{\text{th}}$ local cohomology module of $M$ with respect to the maximal ideal $\fm$.
\end{definition}
\begin{definition} Let $\mathcal F : \ M_0\subset M_1\subset \cdots \subset M_t=M$ be a filtration satisfying the dimension condition. Put $d_i=\dim M_i$. A system of parameters  $\un x=(x_1,$ $\ldots, x_d)$ of $M$ is called a {\it good system of parameters} with respect to $\mathcal F$ if $M_i\cap (x_{d_i+1},\ldots,x_d)M=0$ for $i=0, 1, \ldots, t-1$. A good system of parameters with respect to the dimension filtration is simply called a good system of parameters of $M$.
\end{definition}
The next few results can be implied directly from the definitions or can be found in \cite{cc1}.
\begin{remark}\label{rmq}
 (i) The dimension filtration always exists and it  is unique. In this paper we will always denote the dimension filtration of $M$ by $\D : D_0\subset D_1\subset \cdots \subset D_t=M$.
\item (ii) Let $N\subseteq M$ be a submodule. From the definition of the dimension filtration, there is a $D_i$ such that $N\subseteq D_i$ and $\dim N=\dim D_i$. Consequently, if a filtration $M_0\subset M_1\subset \cdots \subset M_{t^\prime}=M$ satisfies the dimension condition then there exist indices $0\leqslant i_0<i_1<\ldots<i_{t^\prime}$ such that $M_j\subseteq D_{i_j}$ and $\dim M_j=\dim D_{i_j}$. Therefore, a good system of parameters of $M$ is  a good system of parameters with respect to every filtration satisfying the dimension condition.
\item (iii) Let $\mathcal F$ be a filtration satisfying the dimension condition of $M$. Then there always exists on $M$  a good system of parameters with respect to $\mathcal F$. Moreover, if $\un x=(x_1,\ldots ,x_d)$ is a good system of parameters of $M$ with respect to $\mathcal F$, so is $(x_1^{n_1}, \ldots, x_d^{n_d})$ for any integers $n_1,\ldots,n_d>0$.
\item (iv) Let $\un x$ be a good system of parameters. For $\dim D_i<j\leqslant \dim D_{i+1}$, $D_i=0:_Mx_j$. In particular, $0:_Mx_1=H^0_\fm(M)$.
\end{remark}

Let $\mathcal F : \ M_0\subset M_1\subset \cdots \subset M_t=M$ be a filtration satisfying the dimension condition with $d_i=\dim M_i$ and $\un x=(x_1,\ldots,x_d)$ a good system of parameters with respect to $\F$. It is clear that $(x_1,\ldots,x_{d_i})$ is a system of parameters of $M_i$. Therefore the following difference is well defined
$$\ha(\un x)=\ell(M/\un xM)-\sum_{i=0}^te(x_1,\ldots,x_{d_i};M_i),$$
where $e(x_1,\ldots,x_{d_i};M_i)$ is the Serre multiplicity and we set $e(x_1,\ldots,x_{d_0};M_0)=\ell(M_0)$ if  $\dim M_0=0$. Below are some remarkable properties of this number (cf. \cite[Lemma 2.6 and Proposition 2.9]{cc1}).

\begin{lemma}\label{positive} Let $\mathcal F$ be a filtration satisfying the dimension condition and $\un x=(x_1,\ldots,x_d)$ a good system of parameters of $M$. We have
\item(i) $I_{\F, M}(\un x)\geqslant 0$.
\item(ii) Denote $\un x(\un n)=(x_1^{n_1}, \ldots, x_d^{n_d})$ for any $d$-tuple of positive integers  $\un n =(n_1, \ldots, n_d)$ and consider $I_{\F, M}(\un x(\un n))$ as a function in $n_1, \ldots, n_d$, then this function is a non-decreasing function, it means that  $\ha(\un x(\un n))\leqslant \ha(\un x(\un m))$ for all $n_i\leqslant m_i, i=1, \ldots, d$.
\end{lemma}

Concerning  the question of when the function $\ell(M/(x_1^{n_1}, \ldots, x_d^{n_d})M)$ is a  polynomial, the authors in \cite{cc} have introduced a notion of dd-sequences. For the definition we need the notion of d-sequence of Huneke \cite{ch}. A {\it d-sequence} on $M$ is a sequence $(x_1, \ldots, x_s)$ of elements in $\fm$ such that for $i=1, \ldots, s$ and $j\geqslant i$,  $(x_1, \ldots, x_{i-1})M:x_ix_j=(x_1, \ldots, x_{i-1})M:x_j$. 

\begin{definition}\label{ddseq}
A sequence $(x_1, \ldots, x_s)$ of elements in $\fm$  is called a {\it dd-sequence} on $M$ if $(x_1^{n_1}, \ldots, x_i^{n_i})$ is a d-sequence on $M/(x_{i+1}^{n_{i+1}}, \ldots, x_s^{n_s})M$ for all $n_1, \ldots, n_s>0$ and $i=1, \ldots, s$.
\end{definition}
Then dd-sequence is closely related to the notion of good system of parameters by the following lemma.
\begin{lemma}\cite[Lemma 3.5]{cc1}\label{lemma}
Every system of parameters of $M$, which is also a dd-sequence on $M$, is a good system of parameters, and therefore it is a good system of parameters with respect to any filtration $\mathcal F$ satisfying the dimension condition of $M$.
\end{lemma}
We have some characterizations of dd-sequence.
\begin{proposition}\label{tr1}
Let $\un x=(x_1, \ldots, x_d)$ be a system of parameters of $M$. Then the following statements are equivalent:

\item i) $\un x$ is a dd-sequence.

\item ii) For all $0<i\leqslant j\leqslant d$, $n_1, \ldots, n_d>0$, 
$$(x_1^{n_1}, \ldots, \widehat{x_i^{n_i}}, \ldots, \widehat{x_j^{n_j}}, \ldots x_d^{n_d})M:x_i^{n_i}x_j^{n_j}=(x_1^{n_1}, \ldots, \widehat{x_i^{n_i}}, \ldots, \widehat{x_j^{n_j}}, \ldots x_d^{n_d})M:x_j^{n_j}.$$
\item iii) There exist $a_0, a_1, \ldots, a_d \in \mathbb Z$ such that for all $n_1, \ldots, n_d>0$,
$$\ell(M/\un x(\un n)M)=\sum_{i=0}^da_in_1\ldots n_i.$$
In this case, we have $a_i=e(x_1,\ldots, x_i;(x_{i+2},\ldots, x_d)M:x_{i+1}/(x_{i+2},\ldots, x_d)M)$.
\item iv) $\un x$ is a good system of parameters and there exist $b_0, b_1, \ldots, b_{d-1} \in
\mathbb Z$ such that for all $n_1, \ldots, n_d>0$,
$$I_{\D,M}(\un x(\un n))=\sum_{i=0}^{d-1}b_in_1\ldots n_i,$$
where $\D $ is the dimension filtration of $M$
\end{proposition}
\begin{proof}
The implication $(i\Rightarrow ii)$ is proved in \cite[Proposition 3.4]{cc}. For the converse, we need to show that for $0<i\leqslant j<s\leqslant d+1, n_1, \ldots, n_d>0$,
$$(x_1^{n_1}, \ldots, x_{i-1}^{n_{i-1}}, x_s^{n_s}, \ldots, x_d^{n_d})M:x_i^{n_i}x_j^{n_j}=(x_1^{n_1}, \ldots, x_{i-1}^{n_{i-1}}, x_s^{n_s}, \ldots, x_d^{n_d})M:x_j^{n_j},$$
but this is clear by using Krull's Intersection Theorem and the hypothesis.

The equivalence of $(i)$ and $(iii)$ is proved in \cite[Corrollary 3.6]{cc}. By Lemma \ref{lemma}, if $\un x$ is a dd-sequence then it is a good system of parameters. Hence the equivalence of $(iii)$ and $(iv)$ is obvious.
\end{proof}

\begin{lemma} \label{hs1} Let $\D:\ D_0\subset D_1\subset \ldots \subset D_t=M$ be the dimension filtration  and $\un x=(x_1, \ldots, x_d)$ a system of parameters of $M$. Put $d_i=\dim D_i$. Assume that $\un x$ is a dd-sequence on $M$. Then we have $\un xM\cap D_i=(x_1, \ldots, x_{d_i})M\cap D_i$.
\end{lemma}
\begin{proof} We need only to show for any integer $j$, $d_i<j\leqslant d_{i+1}$,  that 
$$\un xM\cap D_i=(x_1,\ldots, x_{j-1}, x_{j+1}, \ldots, x_{d})M\cap D_i.$$ 
Indeed, 
let $a$ be an arbitrary element of $\un xM\cap D_i$. Write $a=x_1a_1+\ldots+x_da_d$. Since $\un x$ is a good system of parameters, ($iv$), $D_i=0:_Mx_{j}$ by Remark \ref{rmq}. Therefore
\begin{multline*}a_{j}\in ((x_1,\ldots, x_{j-1}, x_{j+1}, \ldots, x_d)M+0:_Mx_{j}): x_{j}\\
\subseteq (x_1,\ldots, x_{j-1}, x_{j+1}, \ldots, x_d)M:x_{j}^2=(x_1,\ldots, x_{j-1}, x_{j+1}, \ldots, x_d)M:x_{j},\end{multline*}
and the conclusion follows.
\end{proof}

To end this section, we recall some facts about generalized Cohen-Macaulay modules. For the detailed proof of these results we refer to \cite{cst}. For an $R$-module $M$, we put
$$I(M)=\sup_{\un x}\{\ell(M/\un xM)-e(\un x; M)\},$$
where the supremum is taken over all systems of parameters of $M$. Then $M$ is called a generalized Cohen-Macaulay module if $I(M)<\infty $.  The following characterizations of generalized Cohen-Macaulay modules  are used in this  paper.
\begin{lemma}\label{gcm}
i) If $M$ is a generalized Cohen-Macaulay module, then  $M_\fp$ is Cohen-Macaulay for all $\fp\in \Supp M, \fp\not=\fm$.  Moreover, the converse holds true if $R$ is a factor of a Cohen-Macaulay ring and $M$ is equidimensional.
\vskip 0.2cm
\noindent
ii) The following statements are equivalent:
\item$(1)$ $M$ is a generalized Cohen-Macaulay module.
\item$(2)$ There exist a system of parameters $\un x=(x_1, \ldots, x_d)$ of $M$ and $c\geqslant 0$ such that
$$\ell(M/\un x(\un n)M)=n_1\ldots n_de(\un x;M)+c,$$
for all $n_1, \ldots, n_d>0$. In this case, $c=I(M)$.
\item$(3)$ All the local cohomology modules $H^i_\fm(M)$ are of finite length for $i=0, 1, \ldots, d-1$. 

In particular, if $M$ is a generalized Cohen-Macaulay module then 
$$I(M)=\sum_{i=0}^{d-1}\binom{d-1}{i}\ell(H^i_\fm(M)).$$
\end{lemma}

\section{Sequentially generalized Cohen-Macaulay modules}

First, we recall  the notions of generalized Cohen-Macaulay filtration and  of sequentially generalized Cohen-Macaulay modules, which  were introduced  in \cite{cn}. 
\begin{definition}\label{118}
Let $\F :\ M_0 \subset M_1 \subset \ldots \subset M_t=M$ be a  filtration of submodules of $M$. $\F$ is called a {\it generalized Cohen-Macaulay} filtration if $\F$ satisfies the dimension condition, $\dim M_0\leqslant 0$ and $M_1/M_0, \ldots, M_t/M_{t-1}$ are generalized Cohen-Macaulay modules.

$M$ is called a {\it sequentially generalized Cohen-Macaulay} module if it has a generalized Cohen-Macaulay filtration. 
\end{definition}
By the definition, it is obvious that every generalized Cohen-Macaulay module $M$ is a sequentially generalized Cohen-Macaulay module, where the trivial filtration $0\subset M$ is a generalized Cohen-Macaulay filtration. Suppose that  $M$ is unmixed up to $\fm$-primary, it means that $\dim R/\fp =\dim M$ for all $\fp \in \Ass M\setminus \fm$. Then it is easy to see that $M$ is sequentially generalized Cohen-Macaulay if and only if $M$ is generalized Cohen-Macaulay. Therefore the two-dimensional local domain constructed by Ferrand and Raynaud in \cite{fr} is an example of a  two-dimensional ring which is not a sequentially generalized Cohen-Macaulay module. However, the $\fm$-adic completion of this domain is sequentially generalized Cohen-Macaulay as shown in the following proposition.

\begin{proposition}\label{twodim}
Assume that $R$ is a homomorphic image of a Gorenstein ring and $\dim M=2$. Then $M$ is a sequentially generalized Cohen-Macaulay module.
\end{proposition}
\begin{proof}
Let $N$ be the biggest submodule of $M$ such that $\dim N<2$. Since $ \dim R/\fp=2$ for every $\fp\in \Ass(M/N)$, it is shown by Trung \cite{tr1} that $M/N$ is a generalized Cohen-Macaulay module. If $N$ is of finite length then $M$ has a generalized Cohen-Macaulay filtration $N\subset M$. If $\dim N=1$ then $N$ is generalized Cohen-Macaulay and $M$ has a generalized Cohen-Macaulay filtration $0\subset N\subset M$.
 \end{proof}

The following lemma shows that if $M$ has a generalized Cohen-Macaulay filtration, then it is unique up to $\fm$-primary components and relatively closed to the dimension filtration as follows.
\begin{lemma}\label{cuongnhan}
Let $M$ be a sequentially generalized Cohen-Macaulay module with the dimension filtration $\D : D_0\subset D_1\subset \cdots \subset D_t=M$. Let $\F: M_0\subset M_1\subset\ldots \subset M_{t^\prime}=M$ be a filtration satisfying the dimension condition with $\dim M_1>0$. Then $\F$ is generalized Cohen-Macaulay if and only if $t=t^\prime$ and $\ell(D_i/M_i)<\infty$ for $i=0, 1, \ldots, t-1$. In particular, the dimension filtration of a sequentially Cohen-Macaulay module is always  a generalized Cohen-Macaulay filtration.
\end{lemma}
\begin{proof}
Since $M$ is a sequentially generalized Cohen-Macaulay module, Lemma 4.4 of \cite{cn} shows that the necessary condition holds and $\D$ is a generalized Cohen-Macaulay filtration. We prove the sufficient condition. There are two short exact sequences for each $i=0,1,\ldots, t-1$,
$$0\longrightarrow D_i/M_i\longrightarrow D_{i+1}/M_i\longrightarrow D_{i+1}/D_i\longrightarrow 0,$$
$$0\longrightarrow M_{i+1}/M_i\longrightarrow D_{i+1}/M_i\longrightarrow D_{i+1}/M_{i+1}\longrightarrow 0,$$
where $D_{i+1}/D_i$ is generalized Cohen-Macaulay and $\ell(D_i/M_i)<\infty$. The first exact sequence implies that $D_{i+1}/M_i$ is generalized Cohen-Macaulay. Combining this with the second exact sequence we get that $M_{i+1}/M_i$ is generalized Cohen-Macaulay.
\end{proof}
\begin {remark}\label {dimM1}
Note that without the assumption $\dim M_1>0$  the lemma \ref{cuongnhan} is false. Indeed, if $M_1\not=0$ is of finite length,  both filtrations $0=M_0\subset M_1\subset M_2\subset \ldots \subset M_t=M$ and $M_1\subset M_2\subset \ldots \subset M_t=M$ are generalized Cohen-Macaulay filtrations of lengths $t$ and $t-1$ respectively. For convenience, from now on we only consider generalized Cohen-Macaulay filtrations $M_0\subset M_1\subset\ldots \subset M_t=M$ with $\dim M_1>0$. Then by Lemma \ref{cuongnhan} all generalized Cohen-Macaulay filtrations have the same length which is equal to the length of the dimension filtration. Moreover, Lemma \ref{cuongnhan} enables us to derive many examples of generalized Cohen-Macaulay filtration from a given one. For examples, let $\un x=(x_1, \ldots, x_d)$ be a good system of parameters and $\F: M_0\subset M_1\subset\ldots \subset M_t=M$ a generalized Cohen-Macaulay filtration of $M$.  Then the filtration $0=N_0\subset N_1\subset \ldots \subset N_t=M$ where $N_i=\un x(\un n)M_i, i=1, 2, \ldots, t-1$, is also a generalized Cohen-Macaulay filtration of $M$, and in this example 
$$\ell(M_i/N_i)=\ell(M_i/\un x(\un n)M_i\geqslant n_1\ldots n_{d_i}e(x_1, \ldots, x_{d_i};M_i),$$ 
can be arbitrarily large, where $d_i=\dim M_i$.
\end {remark}
 Note that   a characterization of sequentially generalized Cohen-Macaulay modules by the use  of modules of deficiency was proved in \cite{cn} when $R$ possesses a dualizing complex. In the next, without any restriction on the ground ring, we give a characterization for sequentially generalized Cohen-Macaulay modules by means of local cohomology modules. 
\begin{proposition}\label{cohchar}
$M$ is a sequentially generalized Cohen-Macaulay module if and only if there exists a filtration $\F\ :\ M_0\subset M_1\subset \cdots \subset M_t=M$ satisfying the dimension condition  such that   $\ell(M_0)<\infty$ and $H_\fm^i(M/M_j)$ is of finite length for $j=0, 1, \ldots, t-1$ and $i=0, 1,\ldots, \dim M_{j+1}-1$. Moreover, in this case $\F$ is a generalized Cohen-Macaulay filtration.
\end{proposition}
\begin{proof}
Let $M$ be a sequentially generalized Cohen-Macaulay module with a generalized Cohen-Macaulay filtration $\F:\ M_0\subset M_1\subset \cdots \subset M_t=M$. We prove the necessary condition by induction on the length $t$ of the filtration. The case $t=1$ is proved by Lemma \ref{gcm}, ($ii$). Suppose $t>1$. We observe that $0\subset M_2/M_1\subset \ldots \subset M_{t-1}/M_1\subset M/M_1$ is a generalized Cohen-Macaulay filtration. It follows from the inductive hypothesis that $H_\fm^i(M/M_j)$ is of finite length for  $j=1, \ldots, t-~1$ and $i=0, 1,\ldots, \dim M_{j+1}-1$. It remains to prove that $\ell(H_\fm^i(M))<\infty $ for $i=0, 1, \ldots, \dim M_1-1$. This is clear from the long exact sequence
$$\cdots \longrightarrow H_\fm^i(M_1)\longrightarrow H_\fm^i(M)\longrightarrow H_\fm^i(M/M_1)\longrightarrow \cdots$$
and the fact that $M_1$ is a generalized Cohen-Macaulay module.

For the converse, we consider the long exact sequence
$$\ldots \longrightarrow H_\fm^{i-1}(M/M_j)\longrightarrow H_\fm^i(M_j/M_{j-1})\longrightarrow H_\fm^i(M/M_{j-1})\longrightarrow \ldots$$
Since $H_\fm^{i-1}(M/M_j)$ and $H_\fm^i(M/M_{j-1})$ are of finite length for all $i\leqslant \dim M_j-1$, we have $\ell(H_\fm^i(M_j/M_{j-1}))<\infty$. Hence from Lemma \ref{gcm}, $M_j/M_{j-1}$ is generalized Cohen-Macaulay for $j=1,\ldots, t$.
\end{proof}

Let $\F: M_0\subset M_1\subset \ldots \subset M_t=M$ be a filtration satisfying the dimension condition and  $\un x=(x_1, \ldots, x_d)$ a good system of parameters of $M$ with respect to $\F$. Put $d_i=\dim M_i$. For each $1\leqslant i\leqslant d$, there is $j\in \{0, 1, \ldots, t-1\}$ such that $d_j<i\leqslant d_{j+1}$. We consider the following filtration
$$\F_i:\ (M_0+x_iM)/x_iM\subset \ldots \subset (M_{j-1}+x_iM)/x_iM\subset (M_s+x_iM)/x_iM\subset \ldots \subset M_t/x_iM,$$
where $s=j$ if $d_{j+1}>d_j+1$ and $s=j+1$ if $d_{j+1}=d_j+1$. Then the following lemma is often used in the paper.
\begin{lemma}\label{chia}
Let $M$ be a sequentially generalized Cohen-Macaulay module with a generalized Cohen-Macaulay filtration $\F:\ M_0\subset M_1\subset \ldots \subset M_t=M$. Let $\un x=(x_1, \ldots, x_d)$ be a good system of parameters of $M$ with respect to $\F$.  Then for any  $i\in\{1, 2, \ldots, d\}$, $M/x_iM$ is a sequentially generalized Cohen-Macaulay module with the generalized Cohen-Macaulay filtration $\F_i$ defined as above. 
\end{lemma}
\begin{proof} Let $k$ be a positive integer. If $k<j$, remember the definition of the integer $j$ corresponding to the filtration $\F_i$ we get $d_k<i$, then $M_k\cap x_iM=0$ since $\un x$ is a good system of parameters with respect to $\F$. So $(M_k+x_iM)/(M_{k-1}+x_iM)\simeq M_k/M_{k-1}$ and each quotient module  of the filtration $$(M_0+x_iM)/x_iM\subset \ldots \subset (M_{j-1}+x_iM)/x_iM$$
 is a generalized Cohen-Macaulay module.  Thus, in order to show the generalized Cohen-Macaulay property of the filtration $\F_i$,  it remains to prove that \linebreak $(M_s+x_iM)/(M_{j-1}+x_iM)$ and $(M_k+x_iM)/(M_{k-1}+x_iM)$, $k=s+1, \ldots, t$,  are generalized Cohen-Macaulay.
Let $k\geqslant j$. It is clear that $x_i$ is a parameter element of $M_k$. Let $\D: D_0\subset \ldots \subset D_{t^\prime}=M$ be the dimension filtration of $M$. By Lemma \ref{cuongnhan}, $t^\prime=t$ and $\dim D_k=\dim M_k=d_k$.  Since $\un x$ is a good system of parameters with respect to $\F$ and $\dim(0:_Mx_i)< \dim D_{j+1}$, we have $M_j\subseteq 0:_Mx_i\subseteq D_j$ and $M_k\subseteq 0:_Mx_{d_k+1}\subseteq D_k$. If $x_ia\in M_k\subseteq 0:_Mx_{d_k+1}$ then $a\in 0:_Mx_ix_{d_k+1}\subseteq 0:_Mx_{d_k+1}^2\subseteq D_k$. So $x_iM\cap M_k\subseteq x_iD_k$. We have 
\begin{multline*}\ell\big((x_iM\cap M_k+M_{k-1})/(x_iM_k+M_{k-1})\big)\leqslant \ell\big((x_iD_k+M_{k-1})/(x_iM_k+M_{k-1})\big)\\
\leqslant \ell(x_iD_k/x_iM_k)\leqslant \ell(D_k/M_k)<\infty.\end{multline*}
Thus $(x_iM\cap M_k+M_{k-1})/(x_iM_k+M_{k-1})$ is of finite length. It should be noted that $M_k/(x_iM_k+M_{k-1})\simeq (M_k/M_{k-1})/x_i(M_k/M_{k-1})$ is a generalized Cohen-Macaulay module. Therefore from the short exact sequence
\begin{multline*}0\longrightarrow (x_iM\cap M_k+M_{k-1})/(x_iM_k+M_{k-1})\longrightarrow  M_k/(x_iM_k+M_{k-1})\\
\longrightarrow M_k/(x_iM\cap M_k+M_{k-1})\longrightarrow 0,\end{multline*}
we imply that  $(M_k+x_iM)/(M_{k-1}+x_iM)\simeq M_k/(x_iM\cap M_k+M_{k-1})$ is also generalized Cohen-Macaulay, $k=j, j+1,\ldots, t$. Hence if $s=j$ or equivalently $d_{j+1}>d_j+1$, $\F_i$ is a generalized Cohen-Macaulay filtration. For the case $s=j+1$, that is, $d_{j+1}=d_j+1$, it remains to prove that $(M_{j+1}+x_iM)/(M_{j-1}+x_iM)$ is generalized Cohen-Macaulay. This is immediate by Lemma \ref{gcm}, ($ii$) and the short exact sequence
$$0\longrightarrow M_j/M_{j-1}\longrightarrow M_{j+1}/(x_iM+M_{j-1})\longrightarrow M_{j+1}/(x_iM+M_j)\longrightarrow 0,$$
where $M_j/M_{j-1},\, M_{j+1}/(x_iM+M_j)$ are generalized Cohen-Macaulay modules of dimension $d_j$.
\end{proof}

We say that $M$ is a sequentially Cohen-Macaulay module if each quotient module $D_i/D_{i-1}$ of the dimension filtration $D_0\subset D_1\subset \ldots \subset D_t=M$ of $M$ is a Cohen-Macaulay module, $i=1, \ldots, t$. $M$ is called a locally sequentially Cohen-Macaulay module if for all $\fp\in \Supp M, \fp\not=\fm$, $M_\fp$ is sequentially Cohen-Macaulay. By Lemma \ref{gcm}, a generalized Cohen-Macaulay module is locally Cohen-Macaulay and the converse holds if $R$ is a factor of a Cohen-Macaulay ring and $M$ is equidimensional. There is a similar result for sequentially generalized Cohen-Macaulay modules, however, there is no requirement concerning the equidimensional property of $M$. 
\begin{proposition}
A sequentially generalized Cohen-Macaulay module is locally sequentially Cohen-Macaulay. The converse is true provided $R$ is a factor of a Cohen-Macaulay ring.
\end{proposition}
\begin{proof}
Let $\D:\ D_0\subset D_1\subset \ldots \subset D_t=M$ be the dimension filtration of $M$ and $\fp\in \Supp M, \fp\not=\fm$. Assume that $\Supp M$ is catenary. Let $\fp\in \Supp M$, $\fp\not=\fm$. Using Proposition 2.4 of \cite{sch} we imply that there are $0=i_0 < i_1 < \ldots < i_s\leqslant t$ such that the filtration 
\begin{equation*}\label{loca}
(D_{i_0})_\fp\subset (D_{i_1})_\fp\subset \ldots \subset (D_{i_s})_\fp=M_\fp, \tag{$\star$}
\end{equation*}
is  the dimension filtration of $M_\fp$ and these indices are minimal in the sense $(D_{i_k})_\fp=(D_j)_\fp$, for $i_k\leqslant j<i_{k+1}$, $k=0, 1, \ldots, s$.

Assume that $M$ is a sequentially generalized Cohen-Macaulay module and $\fp \in \Supp M, \fp \not=\fm$. So $\Supp M=\cup_i\Supp D_i/D_{i-1}$ is catenary and $M_\fp$ has the dimension filtration as in (\ref{loca}). Since $D_{i_k}/D_{i_k-1}$ is locally Cohen-Macaulay, $(D_{i_k}/D_{i_{k-1}})_\fp=(D_{i_k}/D_{i_k-1})_\fp$ is Cohen-Macaulay and $M_\fp$ is sequentially Cohen-Macaulay. For the converse, assume in addition that $R$ is a factor of a Cohen-Macaulay ring. Then $R$ is catenary and $M_\fp$ has the dimension filtration as in (\ref{loca}) for $\fp\in \Supp M$, $\fp\not=\fm$. If $\fp\in \Supp D_i/D_{i-1}$ or equivalently, $(D_i)_\fp\not=(D_{i-1})_\fp$, then $i=i_k$ for some $k$ and $(D_i/D_{i-1})_\fp=(D_{i_k}/D_{i_{k-1}})_\fp$ is Cohen-Macaulay since $M_\fp$ is a sequentially Cohen-Macaulay module. Combining this with the fact that $D_i/D_{i-1}$ is equidimentional and $R$ is a factor of a Cohen-Macaulay ring we imply that $D_i/D_{i-1}$ is a generalized Cohen-Macaulay module. So $M$ is a sequentially generalized Cohen-Macaulay module.
\end{proof}

The next result, though its proof is simple,  is the starting point for our study of sequentially generalized Cohen-Macaulay modules in the rest of the paper. 
\begin{theorem}\label{p-stand}
Let $M$ be a sequentially generalized Cohen-Macaulay module with a generalized Cohen-Macaulay filtration $\F: M_0\subset M_1\subset\ldots \subset M_t=M$ and $\un x=(x_1, \ldots, x_d)$ a good system of parameters with respect to $\F$. Then $\ha(\un x(\un n))$ is a constant for all $n_1,\ldots, n_d$  large enough ($n_1,\ldots, n_d\gg 0$ for short).
\end{theorem}
\begin{proof}
Since $\ha(\un x(\un n))$ is non-decreasing by Lemma \ref{positive}, it suffices to prove that  $\ha(\un x(\un n))$ is bounded above by a constant. Put $d_i=\dim M_i$. We have 
\[\begin{aligned}
\ell(M/\un x(\un n)M)&=\ell(M/\un x(\un n)M+M_{t-1})+\ell(\un x(\un n)M+M_{t-1}/\un x(\un n)M)\\
&\leqslant \ell(M/\un x(\un n)M+M_{t-1})+\ell(M_{t-1}/(x_1^{n_1},\ldots, x_{d_{t-1}}^{n_{d_{t-1}}})M_{t-1}).
\end{aligned}\]
Note that $(x_1^{n_1},\ldots, x_{d_i}^{n_{d_i}})$ is a good system of parameter of $M_i$ with respect to the filtration $M_0\subset M_1\subset \ldots \subset M_i$, $i=1, 2, \ldots, t$. By induction on $t$ we have
$$\ell(M/\un x(\un n)M)\leqslant \sum_{i=1}^t\ell(M_i/(x_1^{n_1},\ldots, x_{d_i}^{n_{d_i}})M_i+M_{i-1})+\ell(M_0).$$
Since $M_i/M_{i-1}$ is generalized Cohen-Macaulay, there exists an integer $c\geqslant 0$ such that
\[\begin{aligned}
\ell(M_i/(x_1^{n_1},\ldots,x_{d_i}^{n_{d_i}})M_i+M_{i-1})&\leqslant e(x_1^{n_1},\ldots,x_{d_i}^{n_{d_i}};M_i/M_{i-1})+c\\
&=e(x_1^{n_1},\ldots,x_{d_i}^{n_{d_i}};M_i)+c
\end{aligned}\]
for all $n_1,\ldots, n_d>0$ and $i=1, \ldots, t$. Hence,
$$\ell(M/\un x(\un n)M)\leqslant \sum_{i=1}^te(x_1^{n_1},\ldots,x_{d_i}^{n_{d_i}};M_i)+\ell(M_0)+tc$$
and so $\ha(\un x(\un n))\leqslant tc$ for all $n_1,\ldots, n_d>0$. 
\end{proof}

A consequence of Theorem \ref{p-stand} is the existence of a dd-sequence on a sequentially generalized Cohen-Macaulay module following Remark \ref{rmq}, ($i$) and Proposition \ref{tr1}. Roughly speaking, dd-sequence is another version of p-standard system of parameters defined in \cite{acta}, see also \cite{ta1}, \cite{ta2}. In the case of generalized Cohen-Macaulay module, dd-sequence coincides with the notion of standard system of parameters defined in \cite{tr2}. Standard system of parameters is a powerful tool in studying generalized Cohen-Macaulay modules. p-standard systems of parameters or dd-sequences themselves also have many nice properties and provide a useful tool for studying the structure of non-generalized Cohen-Macaulay modules, see \cite{acta}, \cite{cc}, \cite{cc1}, \cite{ta1}, \cite{ta2}. However, there are examples of modules of which no system of parameters is a dd-sequence. As far as we know, there are only some sufficient conditions for the existence of these systems of parameters, for instance, when the ground ring is a homomorphic image of a Gorenstein ring. The following consequence of Theorem \ref{p-stand} provides another condition.

\begin{corollary} 
Let $M$ be a sequentially generalized Cohen-Macaulay module with a generalized Cohen-Macaulay filtration $\F$ and $\un x=(x_1, \ldots, x_d)$ a good system of parameters of $M$ with respect to $\F$. Then $\un x$ is a dd-sequence if and only if $\ha(\un x(\un n))$ is a constant for all $n_1,\ldots, n_d>0$. In particular, for a sequentially generalized Cohen-Macaulay module $M$ there always exist systems of parameters, which are dd-sequences on $M$. 
\end{corollary}
\begin{proof}
If $\ha(\un x(\un n))$  is a constant then $\un x$ is a dd-sequence by Proposition \ref{tr1}. Vice verse, any dd-sequence is a good system of parameters and $\ha(\un x(\un n))$ is a polynomial in $n_1, \ldots, n_d$. Then $\ha(\un x(\un n))$ must be a constant for all $n_1, \ldots, n_d>0$ by Theorem \ref{p-stand}. Moreover, if $\un x$ is a good system of parameters of $M$, then $\ha(\un x(\un n))$ is non-decreasing and is bounded above by a constant, so it coincides with a constant for $n_1, \ldots, n_d\gg0$. Therefore the existence of a dd-sequence on $M$ follows from the first conclusion and the existence of good system of parameters of $M$.
\end{proof}


\section{The invariant $I_\F(M)$}

Let $M$ be an arbitrary module with a filtration $\F$ satisfying the dimension condition. We put
$$I_\F(M)=\sup_{\un x}{\ha(\un x)},$$
where the supremum is taken over the set of good systems of parameters of $M$ with respect to $\F$. By Theorem \ref{p-stand}, if $\F\ $ is a generalized Cohen-Macaulay filtration and $\un x$ is a good system of parameters with respect to $\F$ then $\ha(\un x(\un n))$ is a constant for all $n_1, \ldots, n_d\gg0$. The aim of this section is to show that this constant  does not depend on the choice of good systems of parameters and  is exactly $I_\F(M)$. Moreover we can compute it by lengths of certain local cohomology modules. 
It should be noticed that when $M$ is a generalized Cohen-Macaulay module and $\F$ is the filtration $0\subset M$, $I_\F(M)$ is exactly the Buchsbaum invariant $I(M)$, which is defined as the supremum of $\ell(M/\un xM)-e(\un x; M)$ taking over all systems of parameters of $M$ (see \cite{sv}). So $I_\F(M)<\infty $ in this case. 
\begin{proposition}\label{bound}
Let $M$ be a sequentially generalized Cohen-Macaulay module and $\F: M_0\subset M_1\subset \ldots \subset M_t=M$ a generalized Cohen-Macaulay filtration of $M$. We have
$$I_\F(M)\leqslant \sum_{i=0}^{t-1}I(M_{i+1}/M_i).$$
In particular, $I_\F(M)<\infty $.
\end{proposition}
\begin{proof}
Let $\un x=(x_1, \ldots, x_d)$ be a good system of parameters of $M$ with respect to $\F$. Put $d_i=\dim M_i$. From the proof of Theorem \ref{p-stand} we obtain 
$$\ell(M/\un xM)\leqslant \sum_{i=0}^{t-1}\ell(M_{i+1}/(x_1, \ldots, x_{d_{i+1}})M_{i+1}+M_i).$$
Hence, 
\[\begin{aligned}
\ha(\un x)&\leqslant \sum_{i=0}^{t-1}\big(\ell(M_{i+1}/(x_1, \ldots, x_{d_{i+1}})M_{i+1}+M_i)-e(x_1, \ldots, x_{d_{i+1}}; M_{i+1})\big)\\
& \leqslant \sum_{i=0}^{t-1}I(M_{i+1}/M_i).
\end{aligned}\]
Taking the supremum of the left hand side over all good systems of parameters with respect to $\F$ we get the result.
\end{proof}
In the next, we will present a computation of $I_\F(M)$ by means of lengths of certain local cohomology modules.  First we need an auxiliary lemma. Recall that a sequence  $(x_1, \ldots, x_s)$ of elements in $\fm$ is said to be a {\it strong d-sequence} on $M$ if $(x_1^{n_1}, \ldots, x_s^{n_s})$ is d-sequence for any $n_1, \ldots, n_s>0$ (see \cite{gy}).
\begin{lemma}\label{vanish}
Let $\un x=(x_1, \ldots, x_d)$ be a system of parameters of $M$ and $N\subset M$ a submodule. Assume that $\un x$ is a strong d-sequence on $M$ and $N\subseteq 0:_Mx_d$. Then we have the following exact sequence for $i<d-1$ and $n\geqslant 3$,
$$0\longrightarrow H^i_\fm(M/N)\longrightarrow H^i_\fm(M/x_d^nM+N)\longrightarrow H^{i+1}_\fm(M/0:x_d)\longrightarrow 0.$$
\end{lemma}
\begin{proof}
Since $\un x$ is a strong d-sequence, the proof of Lemma 2.9 of \cite{acta} implies that $x_jH^i_\fm(M/(x_1, \ldots, x_h)M)=0$ for $j=1, \ldots, d, h+i<j$. So in our case we have $x_dH^i_\fm(M)=0$ for all $i<d$. By then from the long exact sequence
$$\cdots \longrightarrow H^i_\fm(M)\longrightarrow H^i_\fm(M/0:_Mx_d)\longrightarrow H^{i+1}_\fm(0:_Mx_d)\longrightarrow \cdots $$
we obtain $x_d^2H^i_\fm(M/0:_Mx_d)=0$ for all $i<d$. On the other hand, since $0:x_d^n=0:x_d$, we have a commutative diagram
\[\begin{CD}
0@>>>M/0:_Mx_d  @>.x_d^n>> M/N@>>> M/x_d^nM+N@>>>0\\
&&@V.x_d^2VV  @| @VpVV\\
0@>>>M/0:_Mx_d @>.x_d^{n-2}>> M/N@>>> M/x_d^{n-2}M+N@>>>0,
\end{CD}\]
where $p$ is the natural projection. The above diagram derives the following commutative diagram
\[\begin{CD}
\cdots\longrightarrow H_\fm^i(M/0:_Mx_d) @>\psi_i>>   H_\fm^i(M/N) @>>> H_\fm^i(M/x_d^nM+N)\longrightarrow\cdots\\
@V.x_d^2VV @|@VVV\\
\cdots\longrightarrow H_\fm^i(M/0:_Mx_d) @>\varphi_i>> H_\fm^i(M/N) @>>> H_\fm^i(M/x_d^{n-2}M+N)\longrightarrow\cdots,
\end{CD}\]
where $\psi_i, \varphi_i$ are maps derived from the maps $M/0:_Mx_d  \stackrel{.x_d^n}{\longrightarrow} M/N$ and \linebreak $M/0:_Mx_d  \stackrel{.x_d^{n-2}}{\longrightarrow} M/N$ respectively. It implies that $\psi_i=0$ for all $i<d$ since $x_2^2H_\fm^i(M/0:_Mx_d)=0$. So we obtain a short exact sequence for each $i<d-1$,
$$0\longrightarrow H_\fm^i(M/N)\longrightarrow H_\fm^i(M/x_d^nM+N)\longrightarrow H_\fm^{i+1}(M/0:_Mx_d)\longrightarrow 0.$$
\end{proof}

\begin{theorem}\label{invariant}
Let $M$ be a sequentially generalized Cohen-Macaulay module with a generalized Cohen-Macaulay filtration $\F: M_0\subset M_1\subset \cdots \subset M_t=M$. Put $d_i=\dim M_i$. We have
$$I_\F(M)=\ell\big(H^0_\fm(M/M_0)\big)+\sum_{i=0}^{t-1}\sum_{j=1}^{d_{i+1}-1}c_{ij}\ell\big(H^j_\fm(M/M_i)\big),$$
where $c_{ij}=\sum_{k=d_i}^{d_{i+1}-1}\binom{k-1}{j-1}$.
\end{theorem}
         
\begin{proof} 
Let $\un x=(x_1, \ldots, x_d)$ be a good system of parameters of $M$ with respect to $\F$. Since $I_{\F, M}(\un x(\un n))$ is non-decreasing, it suffices to prove that 
$$I_{\F, M}(\un x(\un n))=\ell\big(H^0_\fm(M/M_0)\big)+\sum_{i=0}^{t-1}\sum_{j=1}^{d_{i+1}-1}\sum_{k=d_i}^{d_{i+1}-1}\binom{k-1}{j-1}\ell\big(H^j_\fm(M/M_i)\big),$$
for all $n_1, \ldots, n_d\gg0$. We prove this by induction on the dimension $d$ of $M$. Let $d=1$. Since $x_1$ is a system of parameters of $M$, $x_1^{n_1}M\cap H_\fm^0(M)=0$ for $n_1\gg 0$. So
$$\ell(M/x_1^{n_1}M)=\ell(M/(x_1^{n_1}M+H_\fm^0(M)))+\ell(H_\fm^0(M))=e(x_1^{n_1}, M)+\ell(H_\fm^0(M)).$$
This implies that $\ha(x_1^{n_1})=\ell(H_\fm^0(M))-\ell(M_0)=\ell(H_\fm^0(M/M_0))$, for all $n_1\gg0$. Let $d>1$. By Lemma \ref{chia}, the following filtration is generalized Cohen-Macaulay
$$\F_d:\ (M_0+x_d^{n_d}M)/x_d^{n_d}M\subset \cdots \subset (M_s+x_d^{n_d}M)/x_d^{n_d}M\subset M/x_d^{n_d}M,$$
where $s=t-1$ if $d_{t-1}<d-1$, $s=t-2$ if $d_{t-1}=d-1$. Note that $M_i\cap x_d^{n_d}M=0$, then $(x_d^{n_d}M+M_i)/x_d^{n_d}M\simeq M_i$ for $i=0, 1, \ldots, s$ and $(x_1, \ldots, x_{d-1})$ is a good system of parameters of $M/x_d^{n_d}M$ with respect to $\F_d$. Thus
\[\begin{aligned}
I_{\F_d, M/x_d^{n_d}M}(x_1^{n_1}, \ldots, x_{d-1}^{n_{d-1}})=&\ell(M/\un x(\un n)M)-e(x_1^{n_1}, \ldots, x_{d-1}^{n_{d-1}}; M/x_d^{n_d}M)\\
&-\sum_{i=0}^se(x_1^{n_1}, \ldots, x_{d_i}^{n_{d_i}};M_i).
\end{aligned}\]
On the other hand, since $M/M_{t-1}$ is  generalized Cohen-Macaulay, $(0:_Mx_d^{n_d})/M_{t-1}$ is of finite length. Therefore, if $d_{t-1}=d-1$ then
\[\begin{aligned}
e(x_1^{n_1}, \ldots, x_{d-1}^{n_{d-1}}; M/x_d^{n_d}M)&=e(\un x(\un n);M)+e(x_1^{n_1}, \ldots, x_{d-1}^{n_{d-1}}; 0:_Mx_d^{n_d})\\
&=e(\un x(\un n); M)+e(x_1^{n_1}, \ldots, x_{d-1}^{n_{d-1}}; M_{t-1}).
\end{aligned}\]
Otherwise, if $d_{t-1}<d-1$, $e(x_1^{n_1}, \ldots, x_{d-1}^{n_{d-1}}; M/x_d^{n_d}M)=e(\un x(\un n);M)$. So in both cases, $I_{\F_d, M/x_d^{n_d}M}(x_1^{n_1}, \ldots, x_{d-1}^{n_{d-1}})=I_{\F, M}(\un x(\un n))$ and we have by the inductive hypothesis, 
\[\begin{aligned}
I_{\F, M}(\un x(\un n))
=&\ell\big(H^0_\fm(M/(x_d^{n_d}M+M_0))\big)+\sum_{j=1}^{d-2}\sum_{k=d_s}^{d-2}\binom{k-1}{j-1}\ell\big(H^j_\fm(M/(x_d^{n_d}M+M_s))\big)\\
&+\sum_{i=0}^{s-1}\sum_{j=1}^{d_{i+1}-1}\sum_{k=d_i}^{d_{i+1}-1}\binom{k-1}{j-1}\ell\big(H^j_\fm(M/(x_d^{n_d}M+M_i))\big),
\end{aligned}\]
for all $n_1, \ldots, n_{d-1}\gg0$. By Theorem \ref{p-stand}, $I_{\F, M}(\un x(\un n))$ is a constant for all $n_1, \ldots, n_d\gg0$, hence from Proposition \ref{tr1}, $\un x(\un n)$ is a dd-sequence on $M$ for all $n_1, \ldots, n_d\gg 0$. For $n_d>2$ we can apply Lemma \ref{vanish} to  $M$ and $N=M_j$  to get the following short exact sequence for $i=0, 1, \ldots, t-1$, $j<d-1$,
$$0\longrightarrow H_\fm^j(M/M_i)\longrightarrow H_\fm^j(M/x_d^{n_d}M+M_i)\longrightarrow H_\fm^{j+1}(M/0:_Mx_d^{n_d})\longrightarrow 0.$$
Note that $H_\fm^j(M/x_d^{n_d}M+M_i)\cong H_\fm^{j+1}(M/M_{t-1})$ because $(0:_Mx^{n_d}_d)/M_{t-1}$ is of finite length. Hence $$\ell(H_\fm^j(M/x_d^{n_d}M+M_i))=\ell(H_\fm^j(M/M_i))+\ell(H_\fm^{j+1}(M/M_{t-1}))$$ for $i=0, 1, \ldots, t-1$, $j<d_{i+1}$.
Therefore, for all $n_1, \ldots, n_d\gg 0$,
\[\begin{aligned}
I_{\F, M}(\un x(\un n))=&\ell(H^0_\fm(M/M_0))+\ell(H^1_\fm(M/M_{t-1}))\\
&+\sum_{j=1}^{d-2}\sum_{k=d_s}^{d-2}\binom{k-1}{j-1}\big(\ell(H^j_\fm(M/M_s))+\ell(H^{j+1}_\fm(M/M_{t-1}))\big)\\
&+\sum_{i=0}^{s-1}\sum_{j=1}^{d_{i+1}-1}\sum_{k=d_i}^{d_{i+1}-1}\binom{k-1}{j-1}\big(\ell(H^j_\fm(M/M_i))+\ell(H^{j+1}_\fm(M/M_{t-1}))\big)\\
=&\ell(H^0_\fm(M/M_0))+\sum_{i=0}^{t-1}\sum_{j=1}^{d_{i+1}-1}\sum_{k=d_i}^{d_{i+1}-1}\binom{k-1}{j-1}\ell\big(H^j_\fm(M/M_i)\big).
\end{aligned}\]

\end{proof}
It is proved in \cite[Theorem 1.5]{cc} that $M$ is a sequentially Cohen-Macaulay module if and only if $I_{\D, M}(\un x)=0$ for all good systems of parameters $\un x$ and $\D$ is the dimension filtration of $M$.  In other words, $M$ is a sequentially Cohen-Macaulay module if and only if    $I_\D(M)=0$.  Hence the following characterization of sequentially Cohen-Macaulay modules in terms of local cohomology modules is an immediate consequence of Theorem \ref {invariant}.
\begin{corollary}
Let $\D: D_0\subset D_1\subset \ldots \subset D_t=M$ be the dimension filtration of $M$. $M$ is a sequentially Cohen-Macaulay module if and only if $H^j_\fm(M/D_{i-1})=0$ for all $j<\dim D_i$, $i=1, \ldots, t$.
\end{corollary}

\begin{corollary}
Let $M$ be a sequentially generalized Cohen-Macaulay module, $\F$ a generalized Cohen-Macaulay filtration and $\un x =(x_1, \ldots, x_d)$ a good system of parameters of $M$ with respect to $\F$. Then $\ha(\un x(\un n))\leqslant I_\F(M)$ for all $n_1, \ldots, n_d>0$ and the equality holds for $n_1, \ldots, n_d\gg 0$. In particular, $\un x$ is a dd-sequence on $M$ if and only if $\ha(\un x)=I_\F(M)$. 
\end{corollary}
\begin{corollary}
Let $M$ be a sequentially generalized Cohen-Macaulay module and $\F: M_0\subset M_1\subset \ldots \subset M_t=M$ and  $\F^\prime: N_0\subset N_1\subset \ldots \subset N_t=M$ two generalized Cohen-Macaulay filtrations of $M$. Then
$$I_\F(M)-I_{\F^\prime}(M)=\ell(H^0_\fm(M/M_0))-\ell(H^0_\fm(M/N_0)).$$
\end{corollary}
\begin{proof}
First note that all generalized Cohen-Macaulay filtrations of $M$ have the same length. Let  $D_0\subset D_1\subset \ldots \subset D_t=M$ be the dimension filtration of $M$. From Lemma \ref{cuongnhan}, $D_i/M_i$, $D_i/N_i$ are of finite length for $i=0, 1, \ldots, t$. Hence $H^j_\fm(M/M_i)\simeq H^j_\fm(M/D_i)\simeq H^j_\fm(M/N_i)$ for all $j>0$, $i=0, 1, \ldots, t-1,$ and the conclusion follows from  Theorem \ref{invariant}.
\end{proof}
\begin{corollary}
Let $M$ be a sequentially generalized Cohen-Macaulay module with\linebreak $\mathrm{depth}(M)>0$. Then $I_\F(M)=I_{\F^\prime}(M)$ for two arbitrary generalized Cohen-Macaulay filtrations $\F$,   $\F^\prime$ of $M$. 
\end{corollary}
If $M$ is a generalized Cohen-Macaulay module and $\F: M_0\subset M_1=M$ is a generalized Cohen-Macaulay filtration then it is obvious that $I_\F(M)=\ell(M_0)+I(M/M_0)$. Moreover, we showed in \cite{cc1} that if $\F$ is a Cohen-Macaulay filtration, this means that each $M_{i+1}/M_i$ is Cohen-Macaulay for $i=0, \ldots, t-1$, then $I_\F(M)=\sum_{i=0}^{t-1}I(M_{i+1}/M_i)=0$.  So one might expect that  the inequality in Proposition \ref{bound} becomes an equality in general. Unfortunately, the answer is negative even  $\F$ is the dimension filtration of $M$. We have the following example.
\begin{example}
Let $R=k[[X_1, X_2, X_3, X_4, X_5, X_6]]$ be the ring of all formal power series over a field $k$. We consider the ideals $I=(X_1, X_2, X_3)\cap (X_4, X_5, X_6)$ and $J=(X_2, X_3, X_4, X_5)$. Put $M=R/I\cap J$, then $\dim M=3$. The following filtration is the dimension filtration of $M$,
$$\D:\ 0=D_0\subset D_1\subset D_2=M,$$
where $D_1=I/I\cap J\simeq (I+J)/J\simeq X_1X_6(R/J)$ is  Cohen-Macaulay with $\dim D_1=2$ and $M/D_1=R/I$ is a generalized Cohen-Macaulay module. Therefore $M$ is a sequentially generalized Cohen-Macaulay module. Let $x_1=X_1+X_5, x_2=X_3+X_6, x_3=X_2+X_4$. It could be verified directly that $(x_1, x_2, x_3)$ is a good system of parameters of $M$ and 
$$\ell(M/(x_1^{n_1}, x_2^{n_2}, x_3^{n_3})M)=2n_1n_2n_3+n_1n_2+1,$$
$$\ell((M/D_1)/(x_1^{n_1}, x_2^{n_2}, x_3^{n_3})(M/D_1))=2n_1n_2n_3+2,$$
for all $n_1, n_2, n_3>0$. So $(x_1, x_2, x_3)$ is a dd-sequence on both $M$ and $M/D_1$. Hence, $I_\D(M)=I_{\D, M}(x_1, x_2, x_3)=1$ and $I(M/D_1)=2$. Therefore 
$$I_{\D}(M)=1<0+2 =I(D_1)+I(M/D_1).$$
\end{example}
In the next we will give a necessary and sufficient condition for the inequality mentioned in Lemma \ref{bound} becomes an equality.
\begin{proposition}
Let $M$ be a sequentially generalized Cohen-Macaulay module with a generalized Cohen-Macaulay filtration $\F: M_0\subset M_1\subset \ldots \subset M_t=M$. Then $I_\F(M)=\sum_{i=0}^{t-1}I(M_{i+1}/M_i)$ if and only if we have the following short exact sequences
$$0\longrightarrow H^j_\fm(M_{i+1}/M_i)\longrightarrow H^j_\fm(M/M_i)\longrightarrow H^j_\fm(M/M_{i+1})\longrightarrow 0,$$
for all $i=0, 1, \ldots, t-1$, $j=0, 1, \ldots, \dim M_{i+1}$.
\end{proposition}
\begin{proof}
Denote $d_i=\dim M_i, i=0, 1, \ldots, t$. From the long exact sequence 
$$\ldots \longrightarrow H^j_\fm(M_{i+1}/M_i)\longrightarrow H^j_\fm(M/M_i)\longrightarrow H^j_\fm(M/M_{i+1})\longrightarrow \ldots$$
we imply for $j<d_{i+1}$ that $\ell(H^j_\fm(M_{i+1}/M_i)\geqslant \ell(H^j_\fm(M/M_i)-\ell(H^j_\fm(M/M_{i+1})$. So
\[\begin{aligned}
\sum_{i=0}^{t-1}I(M_{i+1}/M_i)&=\sum_{i=0}^{t-1}\sum_{j=0}^{d_{i+1}-1}\binom{d_{i+1}-1}{j}\ell(H^j_\fm(M_{i+1}/M_i)\\
&\geqslant \sum_{i=0}^{t-1}\sum_{j=0}^{d_{i+1}-1}\binom{d_{i+1}-1}{j}\big(\ell(H^j_\fm(M/M_i)-\ell(H^j_\fm(M/M_{i+1})\big)\\
&=\sum_{i=0}^{t-1}\sum_{j=0}^{d_{i+1}-1}\big(\binom{d_{i+1}-1}{j}-\binom{d_i-1}{j}\big)\ell(H^j_\fm(M/M_i)\\
&=I_\F(M).
\end{aligned}\]
Therefore, $I_\F(M)=\sum_{i=0}^{t-1}I(M_{i+1}/M_i)$ if and only if $\ell(H^j_\fm(M_{i+1}/M_i)= \ell(H^j_\fm(M/M_i)-\ell(H^j_\fm(M/M_{i+1})$ for all $j<d_{i+1}, i=0, 1, \ldots, t-1$. From the above long exact sequence of local cohomology modules again, this is equivalent to the exactness of the following sequences
$$0 \longrightarrow H^j_\fm(M_{i+1}/M_i)\longrightarrow H^j_\fm(M/M_i)\longrightarrow H^j_\fm(M/M_{i+1})\longrightarrow 0,$$
 for all $j<d_{i+1}, i=0, 1, \ldots, t-1$.
\end{proof}
\begin{remark}
In  Theorem \ref{invariant}, the assumption that $\F$ is a generalized Cohen-Macaulay filtration is quite important. For instance, keep all hypothesis in Theorem \ref{invariant}, let $\F^\prime$ be the filtration $0\subset M$. Assume that $\un x$ is a good system of parameters of $M$ such that $\ha(\un x(\un n))$ is constant for all $n_1, \ldots, n_d>0$. We have
$$I_{\F^\prime, M}(\un x(\un n))=\ell(M/\un x(\un n)M)-e(\un x(\un n); M)=\sum_{i=0}^{t-1}n_1 \ldots n_{d_i}e(x_1, \ldots, x_{d_i}; M_i)+I_\F(M).$$
So $I_{\F^\prime}(M)=\infty$ if $t>1$.
\end{remark}


\section{Parametric characterizations}

In the previous sections we have proved the existence of dd-sequence on a sequentially generalized Cohen-Macaulay module and used it to study some properties of these modules. It is shown that $I_\F(M)=\sup_{\un x}\{\ha(\un x)\}$ is finite provided $\F$ is a generalized Cohen-Macaulay filtration of $M$, where the supremum is taken over all good systems of parameters with respect to $\F$.
In this section, we will show that the sequentially generalized Cohen-Macaulayness of $M$ can be characterized by the condition $I_\F(M)<\infty$, where the filtration $\F$ is not necessarily a generalized Cohen-Macaulay filtration. Moreover, we will prove several characterizations of sequentially generalized Cohen-Macaulay property in terms of good systems of parameters. We begin with the following technical lemma.

\begin{lemma}\label{tr}
Suppose that there exist a filtration $\F$ satisfying the dimension condition and a good system of parameters $\un x=(x_1, \ldots, x_d)$ of $M$ with respect to $\F$ such that $\ha(\un x(\un n))$ is a constant for all $n_1, \ldots, n_d>0$. Then $(x_dM:x_i)/(x_dM+0:_Mx_i)$ is of finite length for $i=1, \ldots, d-1$.
\end{lemma}
\begin{proof}
It suffices to prove that $(\un x)(x_dM:x_i)\subseteq x_dM+0:_Mx_i$. Let $\mathcal D : D_0\subset D_1\subset\ldots \subset D_t=M$ be the dimension filtration of $M$ with $d_i=\dim D_i$, $i=0, 1, \ldots, t$. Since $I_{\F, M}(\un x(\un n))=c$ is a constant for all $n_1, \ldots, n_d>0$, we imply by Proposition \ref{tr1} that $\un x$ is a dd-sequence on $M$. Therefore $I_{\D, M}(\un x(\un n))$ is a polynomial in $n_1, \ldots, n_d$. It is clear by Remark \ref{rmq}  ($ii$) and Lemma \ref{positive} that $0\leqslant I_{\mathcal D, M}(\un x(\un n))\leqslant \ha(\un x(\un n))=c$. Thus $I_{\mathcal D, M}(\un x(\un n))$ is also a constant for all $n_1, \ldots, n_d>0$. We prove the lemma by induction on  the length $t$ of the dimension filtration $\D$. Note that the case $i=1$ is trivial since $(x_1,\ldots, x_{d-1})$ is a d-sequence on $M/x_dM$. If $t=1$ then $I_{\D, M}(\un x(\un n))=\ell(M/\un x(\un n)M)-e(\un x(\un n);M)-\ell(D_0)$ is constant for all $n_1, \ldots, n_d>0$. Thus $(x_1, \ldots, x_d)$ is a dd-sequence in any order. Hence $x_dM:x_i=x_dM:x_j$ for all $i, j<d$ and $(\un
x)(x_dM:x_i)\subseteq x_dM\subseteq x_dM+0:_Mx_i$.

For each $t>1$ we prove the assertion by induction on $d_1$. Let $d_1=1$.   Since $D_1=0:_Mx_2$ and $(x_1, x_2)$ is a strong d-sequence, it follows by Lemma \ref{hs1} that $$D_1\cap \un x(\un n)M=D_1\cap x_1^{n_1}M=x_1^{n_1}D_1.$$ Hence 
$$\ell(M/\un x(\un n)M+D_1) =\ell(M/\un x(\un n)M)-\ell(D_1/x_1^{n_1}D_1)=\sum_{i=2}^te(x_1^{n_1},\ldots, x_{d_i}^{n_{d_i}}; D_i)+c.$$
Note that $M/D_1$ has the dimension filtration $$\D^\prime : 0\subset D_2/D_1\subset \cdots \subset D_t/D_1=M/D_1$$ 
and $e(x_1^{n_1},\ldots, x_{d_i}^{n_{d_i}}; D_i)=e(x_1^{n_1},\ldots, x_{d_i}^{n_{d_i}}; D_i/D_1),$ $i>1$. Thus $I_{\D^\prime, M/D_1}(\un x(\un n))=c$. Applying the inductive hypothesis to $\D^\prime$ we obtain
$$(x_1, \ldots, x_{d-1})[x_d(M/D_1):x_i]\subseteq x_d(M/D_1)+(0:x_i)_{M/D_1}$$ for all $1<i<d$. Thus 
$$(x_1, \ldots, x_{d-1})[(x_dM+D_1):x_i]\subseteq x_dM+D_1:x_i= x_dM+0:_Mx_i,$$
since $\un x$ is a d-sequence on $M$. So $(x_1, \ldots, x_d)(x_dM:x_i)\subseteq x_dM+0:_Mx_i$.

Assume $d_1>1$. It is easy to check that the following filtration of $M/x_1^{n_1}M$  satisfies  the dimension condition 
$$\D_1 : (x_1^{n_1}M+D_0)/x_1^{n_1}M\subset \ldots \subset (x_1^{n_1}M+D_{t-1})/x_1^{n_1}M\subset M/x_1^{n_1}M,$$
where $(x_1^{n_1}M+D_i)/x_1^{n_1}M\simeq D_i/D_i\cap x_1^{n_1}M= D_i/x_1^{n_1}D_i$. Hence
\[\begin{aligned}
e(x_2^{n_2},\ldots,x_{d_i}^{n_{d_i}}; (x_1^{n_1}M+D_i)/x_1^{n_1}M)
=&e(x_2^{n_2},\ldots,x_{d_i}^{n_{d_i}}; D_i/x_1^{n_1}D_i)\\
=&e(x_1^{n_1}, \ldots, x_{d_i}^{n_{d_i}}; D_i).
\end{aligned}\]
Therefore $I_{\D_1, M/x_1^{n_1}M}(x_2^{n_2},\ldots, x_d^{n_d})=I_{\D, M}(\un x(\un n))=c$ for all $n_1,\ldots, n_d$. Note that $\dim(x_1^{n_1}M+D_1)/x_1^{n_1}M=\dim D_1/x_1^{n_1}D_1=d_1-1$. Using the inductive hypothesis we obtain
$$(x_2, \ldots, x_d)[x_d(M/x_1^{n_1}M):x_i]\subseteq x_d(M/x_1^{n_1}M)+(0:x_i)_{M/x_1^{n_1}M}.$$
In other words, $(x_2, \ldots, x_d)[(x_d, x_1^{n_1})M:x_i]+x_1^{n_1}M\subseteq x_dM+x_1^{n_1}M:x_i$. By Krull's Intersection Theorem we have
\[\begin{aligned}
(x_2, \ldots, x_d)(x_dM:x_i)&\subseteq \bigcap_{n_1}\big((x_2, \ldots, x_d)[(x_d, x_1^{n_1})M:x_i]+x_1^{n_1}M\big)\\
&\subseteq \bigcap_{n_1}(x_dM+x_1^{n_1}M:x_i)=x_dM+0:_Mx_i.
\end{aligned}\]
We also have $I_{\D_1, M/x_1^{n_1}M}(x_2^{n_2},x_1^{n_1},x_3^{n_3},\ldots, x_d^{n_d})=c$ since $d_1\geqslant 2$ . Applying the same method to the sequence $(x_2, x_1, x_3, \ldots, x_d)$ we get $(x_1, x_3, \ldots, x_d)(x_dM:x_i)\subseteq x_dM+0:_Mx_i$. So $(\un x)(x_dM:x_i)\subseteq x_dM+0:_Mx_i$ as required.
\end{proof}

\begin{theorem}\label{main}
Let $M$ be a finitely generated $R$-module of dimension $d$. The following statements are equivalent:

i) $M$ is a sequentially  generalized Cohen-Macaulay module.

ii) There exists a filtration $\F$ of submodules of $M$ satisfying the dimension condition such that $I_\F(M)<\infty$.

iii) There exists a filtration $\F$ of submodules of $M$ satisfying the dimension condition and a good system of parameters $\un x=(x_1, \ldots, x_d)$ of $M$ with respect to $\F$ such that $\ha(\un x(\un n))$ is a constant for all $n_1, \ldots, n_d>0$.
\end{theorem}
\begin{proof} $(i)\Rightarrow (ii)$ is the content of Theorem \ref{invariant}.\vspace{0.15cm}

$(ii)\Rightarrow (iii)$ is straightforward since $\ha(\un x(\un n))$ is non-decreasing. \vspace{.15cm}

$(iii)\Rightarrow (i)$: Let $\mathcal D: \ D_0\subset \ldots \subset D_t=M$ be the dimension filtration of $M$. By the same argument as in the proof of Lemma \ref{tr} we get that $I_{\mathcal D, M}(\un x(\un n))$ is a constant for all $n_1,\ldots, n_d>0$. 
Now, we argue  the statement by induction on $d$. The case $d=1$ is trivial since $M$ is a generalized Cohen-Macaulay module.
Assume that $d>1$. Consider the following filtration of $M/x_dM$
$$\D_d : (x_dM+D_0)/x_dM\subset \cdots \subset  (x_dM+D_s)/x_dM\subset M/x_dM,$$
where $s=t-1$ if $d_{t-1}<d-1$ and $s=t-2$ if $d_{t-1}=d-1$. Since $\un x$ is a good system of parameters of $M$, $D_i\cap x_dM=0$ and $(D_i+x_dM)/x_dM\simeq D_i$ for all $i=0, 1, \ldots, t-1$. So $\dim(D_i+x_dM)/x_dM=d_i$ and $\D_d$ satisfies  the dimension condition. Since $(x_1, \ldots, x_{d-1})$ is a dd-sequence on $M/x_dM$, it is a good
system of parameters of $M/x_dM$ by Lemma \ref{lemma}. It is not difficult to verify that  $I_{\D_d, M/x_dM}(x_1^{n_1},\ldots, x_{d-1}^{n_{d-1}})= I_{\D, M}(\un x(\un n))=c$ for a non-negative constant $c$ and all $n_1,\ldots, n_{d-1}>0, n_d=1$. Therefore, from the inductive hypothesis $M/x_dM$ is a sequentially generalized Cohen-Macaulay module.  By Remark
\ref{rmq}, ($iv$), there is a sequence of non-negative integers $l_0<l_1<\ldots<l_r$ such that the filtration 
$$\overline{\mathcal D}: \overline D_0=(0:x_{l_0+1})_{M/x_dM}\subset \ldots \subset
\overline D_r=(0:x_{l_r+1})_{M/x_dM}\subset M/x_dM$$
is the dimension
filtration of $M/x_dM$ with
 $\dim \overline D_k=l_k$. By Remark \ref{rmq}, (ii), for each $0\leqslant i\leqslant s$ there is an $k$ such that $l_k=d_i$ and
$(D_i+x_dM)/x_dM\subseteq \overline D_k.$
Hence,
$$e(x_1, \ldots, x_{l_k};\overline D_k)\geqslant e(x_1, \ldots,
x_{d_i};(D_i+x_dM)/x_dM).$$
On the other hand, since
$$0\leqslant I_{\overline{\mathcal D}, M/x_dM}(x_1^{n_1}, \ldots,
x_{d-1}^{n_{d-1}})\leqslant I_{\mathcal D_d, M/x_dM}(x_1^{n_1}, \ldots,
x_{d-1}^{n_{d-1}})=c$$ for all $n_1, \ldots, n_{d-1}>0$, it follows that
\begin{equation*}\label{thang}
\sum_{k=0}^rn_1\ldots n_{l_k}e(x_1, \ldots,
x_{l_k};\overline D_k)\leqslant c+\sum_{i=0}^sn_1\ldots
n_{d_i}e(x_1, \ldots, x_{d_i};(D_i+x_dM)/x_dM).\end{equation*}

Therefore we obtain that  $r=s$ and $i_k=d_k$ for $k=0, 1, \ldots, s$.
It should be noted that $$\overline D_i/(x_dM+D_i/x_dM)\simeq (x_dM:x_{d_i+1})/(x_dM+D_i)= (x_dM:x_{d_i+1})/(x_dM+0:_Mx_{d_i+1})$$ is of finite length by Lemma \ref{tr}, and so $\D_d$ is a generalized Cohen-Macaulay filtration by Proposition \ref{cuongnhan}. Thus each quotient $D_i/D_{i-1}$ for $i=1,\ldots , s$  and $M/x_dM+D_s$ are generalized Cohen-Macaulay modules. Now, replace $x_d$ by $x_d^3$. We have to consider two cases.

\noindent {\it Case 1}. $d_{t-1}<d-1$, then $s=t-1$ and it remains to prove that $M/D_{t-1}$ is a generalized Cohen-Macaulay module. Applying Lemma \ref{vanish} to the module $M$ with $N=D_{t-1}$ and the dd-sequence  $\un x$, we have the following short exact sequence for $i<d$,
$$0\longrightarrow H^{i-1}_\fm(M/D_{t-1})\longrightarrow H^{i-1}_\fm(M/x_d^3M+D_{t-1})\longrightarrow H^i_\fm(M/D_{t-1})\longrightarrow 0.$$
We have just proved that $M/x_d^3M+D_{t-1}$ is a generalized Cohen-Macaulay module. Therefore $\ell(H^i_\fm(M/D_{t-1}))\leqslant \ell(H^{i-1}_\fm(M/x_d^3M+D_{t-1}))<\infty$, $i=1, 2, \ldots, d-1$, and $M/D_{t-1}$ is a generalized Cohen-Macaulay module.

\noindent {\it Case 2}. $d_{t-1}=d-1$, then $s=t-2$. We need to prove that $M/D_{t-1}$ and $D_{t-1}/D_{t-2}$ are generalized Cohen-Macaulay. Using Lemma \ref{vanish} for $M$ and $N=D_{t-2}$ we have a short exact sequence
$$0\longrightarrow H^{i-1}_\fm(M/D_{t-2})\longrightarrow H^{i-1}_\fm(M/x_d^3M+D_{t-2})\longrightarrow H^i_\fm(M/D_{t-1})\longrightarrow 0$$
for all $i<d$. Since $M/x_d^3M+D_{t-2}$ is generalized Cohen-Macaulay, $\ell(H^i_\fm(M/D_{t-1}))\leqslant \ell(H^{i-1}_\fm(M/x_d^3M+D_{t-2}))<\infty$. Therefore $M/D_{t-1}$ is generalized Cohen-Macaulay. It should be noted that $M_{t-1}\cap x_dM=0$. We have a short exact sequence
$$0\longrightarrow D_{t-1}/D_{t-2} \longrightarrow M/x_dM+D_{t-2}\longrightarrow M/x_dM+D_{t-1}\longrightarrow 0.$$
Since $M/x_dM+D_{t-1}$ and $M/x_dM+D_{t-2}$ are both generalized Cohen-Macaulay of dimension $d-1$,  so is $D_{t-1}/D_{t-2}$, and the proof of Theorem \ref{main} is complete.
\end{proof}
It is known that $M$ is a generalized Cohen-Macaulay module if and only if for every system of parameters $\un x=(x_1, \ldots, x_d)$, $(x_1^{n_1},\ldots, x_d^{n_d})$ is a d-sequence for all $n_1, \ldots, n_d\gg 0$. This result raises a nature question: whether $M$ is a sequentially Cohen-Macaulay module if there is a filtration $\F$ such that for every good system of parameters $\un x$ with respect to $\F$, $(x_1^{n_1}, \ldots, x_d^{n_d})$ is a dd-sequence on $M$ for all $n_1,\ldots,n_d\gg 0$. Unfortunately the answer is negative as in the following example.
\begin{example}
Let $S=k[[x,y,z,t,w]]$ be the ring of formal power series with coefficients in a field $k$. Put $R=S/(yt,yw,zt,zw)$. Then $\dim R=3$ and the non-Cohen-Macaulay locus of $R$ is
$$\NCM(R)=\{\fp\in \Spec R : R_\fp \text{ is not Cohen-Macaulay}\}=V(y,z,t,w).$$
Put $M_1=R/(y,z,t,w)$ and $M=M_1\oplus R$. $M$ has the dimension filtration $\mathcal D: 0\subset M_1\subset M$. $\mathcal D$ is not a generalized Cohen-Macaulay filtration since $M/M_1\simeq R$ is not a generalized Cohen-Macaulay ring ($\dim \NCM(R)=1$). Let $\un x=(x_1,x_2,x_3)$ be a good system of parameters of $M$. Then $x_2, x_3\in \Ann M_1=(y,x,t,w)=\Rad(\fa(R))$ where $\fa(R)=\Ann H_\fm^0(R)\Ann H_\fm^1(R)\Ann H_\fm^2(R)$. By \cite[Corollary 3.9]{acta}, $(x_1^{n_1},x_2^{n_2},x_3^{n_3})$ is a dd-sequence for all $n_1,n_2,n_3\gg 0$.
\end{example}
In the next example, we want to clarify that the filtration $\F$ mentioned in Theorem \ref{main} do not need  to be a generalized Cohen-Macaulay filtration.
\begin{example}
Let $R=k[[x, y, z, w]]$ be the ring of formal power series over a field
$k$. We put $M=R/(xy, xz)$ and $M_1=(xy, xz, xw)/(xy, xz)$. Then $\dim
M=3$, $\dim M_1=2$ and the filtration $\F: 0\subset M_1\subset M$
satisfies the dimension condition. Note that $M/M_1\simeq R/(xy, xz,
xw)=R/(x)\cap(y, z, w)$ is not a generalized Cohen-Macaulay module, thus
$\F$ is not a generalized Cohen-Macaulay filtration. On the other hand, it
is easy  to verify that $(w, x+y, z)$ is a good system of parameters of $M$
with respect to $\F$ and
$$\ell(M/(w^l, (x+y)^m, z^n)=lmn+lm=lmne(w, x+y, z; M)+lme(w, x+y; M_1).$$
In other words, $I_{\F, M}(w^l, (x+y)^m, z^n)=0$ for all $l, m, n>0$. Thus $M$ is a sequentially generalized Cohen-Macaulay module by Theorem \ref{main}.

More general, let $M$ be a sequentially generalized Cohen-Macaulay module with the dimension filtration $\D:\ D_0\subset D_1\subset \ldots \subset D_t=M$. By Lemma \ref{cuongnhan}, $\D$ is a generalized Cohen-Macaulay filtration. Let $\un x=(x_1, \ldots, x_d)$ be a good system of parameters of $M$. We consider the following filtration $\F:\ 0=M_0\subset M_1\subset \ldots \subset M_t=M$ where $M_i=x_1D_i$ for all $0<i<t$. Put $d_i=\dim D_i$. Since $(x_1, \ldots, x_{d_i})$ is a system of parameters of $D_i$, $\dim D_i/x_1D_i=d_i-1$ for all $i>0$. So by Lemma \ref{cuongnhan}, $\F$ is not a generalized Cohen-Macaulay filtration if $t>1$. On the other hand, $e(x_1, \ldots, x_{d_i}; M_i)=e(x_1, \ldots, x_{d_i}; D_i)$ for all $i>0$ and 
$$I_{\F, M}(\un x(\un n))=I_{\D, M}(\un x(\un n))+\ell(D_0)$$
which is bounded above by a constant for all $n_1, \ldots, n_d>0$.
\end{example}

The following theorem gives a finite criterion for the sequentially generalized Cohen-Macaulay property.
\begin{theorem}\label{finite}
A finitely generated $R$-module $M$ is a sequentially generalized Cohen-Macaulay module if and only if there exist a filtration $\F : M_0\subset M_1\subset \cdots \subset M_t=M$ satisfying the dimension condition and a good system of parameters $\un x=(x_1, \ldots, x_d)$ with respect to $\F$ such that $\ha(x_1,\ldots, x_d)=\ha(x_1^2,\ldots, x_d^2)$.
\end{theorem}
\begin{proof}
Put $\ha(x_1,\ldots, x_d)=c$. By Theorem \ref{main}  it suffices to prove that $\ha(\un x(\un n))=c$ for all $n_1, \ldots, n_d>0$. The proof is established by induction on the dimension of $M$. The case $d=1$ is immediate because $M$ is a generalized Cohen-Macaulay module. Assume $d>1$. Since the function $\ha(\un x(\un n))$ is non-decreasing, we have $\ha(\un x(\un n))=c$ for all  $1\leqslant n_1,\ldots, n_d\leqslant 2$. We first prove that $\ha(\un x(\un n))$ does not depend on $n_d$ for a fixed $(d-1)$-tuple $(n_1,\ldots, n_{d-1})$ with $1\leqslant n_1,\ldots, n_{d-1}\leqslant 2$. We have 
$$\ell(M/\un x(\un n)M)-e(\un x(\un n);M)=\ha(\un x(\un n))+\sum_{i=0}^{t-1}e(x_1^{n_1},\ldots, x_{d_i}^{n_{d_i}};M_i),$$
 which is independent of $n_d$ for $n_d\in\{1,2\}$ by the hypothesis. Applying Corollary 4.3 of \cite{ab} to $M$ and the system of parameters $(x_d^{n_d}, x_1^{n_1},\ldots, x_{d-1}^{n_{d-1}})$ we have
\[\begin{aligned}
\ell(M/\un x(\un n)M)-e(\un x(\un n);M)=&\sum_{i=1}^{d-1}n_de\big(x_d,x_1^{n_1},\ldots,x_{i-1}^{n_{i-1}}; (0:x_i^{n_i})_{M/(x_{i+1}^{n_{i+1}}, \ldots, x_d^{n_d})M}\big)\\&+\ell\big((0:x_d^{n_d})_{M/(x_1^{n_1},\ldots,x_{d-1}^{n_{d-1}})M}\big),
\end{aligned}\]
which is non-decreasing in $n_d$. Let $n_d$ vary in $\{1, 2\}$, we get
$$e\big(x_d,x_1^{n_1},\ldots,x_{i-1}^{n_{i-1}}; (0:x_i^{n_i})_{M/(x_{i+1}^{n_{i+1}}, \ldots, x_d^{n_d})M}\big)=0 $$ for all $0< i< d$ and
$$(0:x_d^2)_{M/(x_1^{n_1},\ldots,x_{d-1}^{n_{d-1}})M}=(0:x_d)_{M/(x_1^{n_1},\ldots,x_{d-1}^{n_{d-1}})M}.$$
The last equality implies 
$$(0:x_d^{n_d})_{M/(x_1^{n_1},\ldots,x_{d-1}^{n_{d-1}})M}=(0:x_d)_{M/(x_1^{n_1},\ldots,x_{d-1}^{n_{d-1}})M}$$
for all $n_d>0$. So we have
\[\begin{aligned}\ell(M/\un x(\un n)M)-e(\un x(\un n);M)=&\ell\big((0:x_d)_{M/(x_1^{n_1},\ldots,x_{d-1}^{n_{d-1}})M}\big)\\=&\ell(M/(x_1^{n_1},\ldots, x_{d-1}^{n_{d-1}}, x_d)M)-e(x_1^{n_1},\ldots, x_{d-1}^{n_{d-1}}, x_d;M)\\
=&c+\sum_{i=0}^{t-1}e(x_1^{n_1},\ldots, x_{d_i}^{n_{d_i}}; M_i),
\end{aligned}\]
and $\ha(\un x(\un n))=c$ for all $1\leqslant n_1, \ldots, n_{d-1}\leqslant 2$ and all $n_d>0$.

Put $\un x^\prime(\un n)=(x_1^{n_1},\ldots, x_{d-1}^{n_{d-1}})$. It is not difficult to verify that the following filtration satisfies the dimension condition
$$\F_d : (M_0+x_d^{n_d}M)/x_d^{n_d}M\subset \cdots \subset (M_s+x_d^{n_d}M)/x_d^{n_d}M\subset M/x_d^{n_d}M,$$
where $s=t-2$ if $d_{t-1}<d-1$ and $s=t-1$ if $d_{t-1}=d-1$, and 
$$c=\ha(\un x(\un n))=I_{\F_d,M/x_d^{n_d}M}(\un x^\prime(\un n))+e(\un x^\prime(\un n);0:_Mx_d^{n_d}/M_{t-1})$$
for all $n_1, \ldots, n_{d-1}\in \{1, 2\}, n_d>0$. Note that each term in the right  of this equality is non-decreasing in $n_1, \ldots, n_{d-1}$, thus $e(\un x^\prime(\un n);0:_Mx_d^{n_d}/M_{t-1})=0$ and $I_{\F_d,M/x_d^{n_d}M}(\un x^\prime(\un n))=c$ for all $n_1, \ldots, n_{d-1}\in \{1, 2\}$. So by the inductive hypothesis $I_{\F_d,M/x_d^{n_d}M}(\un x^\prime(\un n))=c$ for all $n_1, \ldots, n_{d-1}>0$. Therefore $\ha(\un x(\un n))=I_{\F_d,M/x_d^{n_d}M}(\un x^\prime(\un n))=c$ for all $n_1, \ldots, n_d>0$.
\end{proof}
In many cases, it is not easy to verify that a system of parameters is a dd-sequence or not. To do this we usually use one of the equivalent conditions stated in Lemma \ref{tr1}. Theorem \ref{finite} and its proof provide another finite criterion for examining whether a system of parameters is a dd-sequence on a sequentially generalized Cohen-Macaulay module.
\begin{corollary}
Let $M$ be a sequentially generalized Cohen-Macaulay module. A system of parameters $\un x=(x_1, \ldots, x_d)$ of $M$ is a dd-sequence on $M$ if and only if there exists a filtration $\F$ satisfying the dimension condition such that $\un x$ is good with respect to $\F$ and $\ha(\un x)=\ha(x_1^2, \ldots, x_d^2)$.
\end{corollary}


\section{Hilbert-Samuel function}
It has been shown in Section 3 that any sequentially generalized Cohen-Macaulay module $M$ admits a dd-sequence, i. e., a good system of parameters $\un x=(x_1, \ldots, x_d)$ such that $I_{\D, M}(\un x(\un n))$ is a constant for all $n_1, \ldots, n_d>0$, where $\D$ is the dimension filtration of $M$. Denote $\fq=(x_1, \ldots, x_d)R$. The aim of this section is to study the Hilbert-Samuel function of $M$ with respect to   $\fq$. We show that when $M$ is a sequentially generalized Cohen-Macaulay module and $\un x$ is a dd-sequence on $M$, this function coincides with the Hilbert-Samuel polynomial. Moreover, the coefficients of this polynomial might be expressed in terms of lengths of certain local cohomology modules. 
\begin{lemma}\label{hs2}
Let $M$ be a sequentially generalized Cohen-Macaulay module with a generalized Cohen-Macaulay filtration $\mathcal F: M_0\subset M_1\subset \ldots \subset M_t=M$. Let $\un x=(x_1, \ldots, x_d)$ be a good system of parameters of $M$ with respect to $\F$, which is a dd-sequence. Then we have the following short exact sequences
$$0\longrightarrow H_\fm^i(M)\longrightarrow H_\fm^i(M/x_1M)\longrightarrow H_\fm^{i+1}(M) \longrightarrow 0,$$
for $0\leqslant i\leqslant \dim M_1-2$.
\end{lemma}
\begin{proof}
Put $d_1=\dim M_1$. It is obvious that $(x_2, x_3, \ldots, x_{d_1}, x_1, x_{d_1+1}, \ldots, x_d)$ is also a good system of parameters of $M$ with respect to $\F$, and 
$$I_{\D,M}(x_2^{n_2}, x_3^{n_3}, \ldots, x_{d_1}^{n_{d_1}}, x_1^{n_1}, x_{d_1+1}^{n_{d_1+1}}, \ldots, x_d^{n_d})=I_{\D, M}(\un x(\un n))$$
is a constant for all $n_1, \ldots, n_d>0$ by Theorem \ref{p-stand}. It  follows from Proposition \ref{tr1} that $(x_2, x_3, \ldots, x_{d_1},$ $x_1, x_{d_1+1}, \ldots, x_d)$ is a dd-sequence on $M$, and hence it is a strong d-sequence on $M$. Then $x_1H^i_\fm(M)=0$ for all $i<d_1$ (see \cite[Lemma 2.9]{acta}) and $0:_Mx_1=0:_M(\un x)R\subseteq H^0_\fm(M)$ is of finite length. Therefore from the long exact sequence of local cohomology module
\begin{multline*}
0\longrightarrow H^0_\fm(M)\longrightarrow H^0_\fm(M/x_1M)\longrightarrow H^1_\fm(M)\stackrel{.x_1}{\longrightarrow} H^1_\fm(M)\longrightarrow \ldots \\
\longrightarrow H^i(M)\stackrel{.x_1}{\longrightarrow}H^i_\fm(M)\longrightarrow H^i_\fm(M/x_1M)\longrightarrow H^{i+1}_\fm(M)\stackrel{.x_1}{\longrightarrow} \ldots
\end{multline*}
we obtain the short exact sequences
$$0\longrightarrow H_\fm^i(M)\longrightarrow H_\fm^i(M/x_1M)\longrightarrow H_\fm^{i+1}(M)\longrightarrow 0,$$
for $0\leqslant i\leqslant d_1-2$.
\end{proof}

\begin{theorem}\label{hs6}
Let $M$ be a sequentially generalized Cohen-Macaulay module with a generalized Cohen-Macaulay filtration $\F: M_0\subset M_1\subset \ldots \subset M_t=M$ and $\un x=(x_1, \ldots, x_d)$ a system of parameters of $M$. Assume that $\un x$ is a dd-sequence on $M$. Put $d_i=\dim M_i$ and $\fq=(x_1, \ldots, x_d)R$. Then for all $n\geqslant 0$ we have
\begin{equation}\label{1}\ell(M/\fq^{n+1}M)=\sum_{i=0}^d\binom{n+i}{i}e_{d-i}(\fq; M),\tag{1}\end{equation}
where $e_d(\fq;M)=\ell(H^0_\fm(M))$,
\begin{equation}\label{2}
e_{d-d_k}(\fq;M)=e(x_1, \ldots, x_{d_k}; M_k)+\sum_{j=1}^{d_k}\binom{d_k-1}{j-1} \ell(H_\fm^j(M/M_k)), \tag{2}
\end{equation}
for $k=1, \ldots, t$,
and 
\begin{equation}\label{3}
e_{d-i}(\fq;M)=\sum_{j=1}^i\binom{i-1}{j-1}\ell(H_\fm^j(M/M_k))\tag{3}
\end{equation}
for $d_k<i<d_{k+1}$.
\end{theorem}
\begin{proof}
Since $\un x$ is a d-sequence, it was shown by Trung in \cite[Theorem 4.1]{tr11} that the Hilbert-Samuel function $\ell(M/\fq^{n+1}M)$ admits the expression (\ref{1}) where by a slight modification, $e_d(\fq;M)=\ell(H^0_\fm(M))$ and
$$e_{d-i}(\fq;M)=\ell\big(H^0_\fm(M/(x_1, \ldots, x_i)M)\big)-\ell\big(H^0_\fm(M/(x_1, \ldots, x_{i-1})M)\big),\ i>0.$$
We argue (\ref{2}), (\ref{3}) by induction on $d$. The case $d=1$ is trivial. Let $d>1$. Firstly assume that $d_1>1$. We have
$$e_{d-1}(\fq;M)=\ell(H^0_\fm(M/x_1M))-\ell(H^0_\fm(M)).$$
From Lemma \ref{hs2} there is a short exact sequence
$$0\longrightarrow H_\fm^0(M)\longrightarrow H_\fm^0(M/x_1M)\longrightarrow H_\fm^1(M)\longrightarrow 0.$$
Then $e_{d-1}(\fq;M)=\ell(H^1_\fm(M))$. By Lemma \ref{chia} the following filtration 
$$\F_1: (M_0+x_1M)/x_1M\subset (M_1+x_1M)/x_1M\subset \ldots \subset M/x_1M$$
is a generalized Cohen-Macaulay filtration of $M/x_1M$.
Hence from the inductive hypothesis we get the following equality for $k=1, \ldots, t$,
\begin{multline*}
e_{d-d_k}(\fq;M)=e_{d-d_k}(x_2, \ldots, x_d; M/x_1M)=e(x_2, \ldots, x_{d_k}; (M_k+x_1M)/x_1M)\\
+\sum_{j=1}^{d_k-1}\binom{d_k-2}{j-1}\ell(H^j_\fm(M/(x_1M+M_k))),\end{multline*}
and 
$$e_{d-i}(\fq;M)=e_{d-i}(x_2, \ldots, x_d; M/x_1M)=\sum_{j=1}^{i-1}\binom{i-2}{j-1}\ell(H^j_\fm(M/(x_1M+M_k)))$$
for $d_k<i<d_{k+1}$.
Using Lemma \ref{hs2} again we obtain
$$e_{d-d_k}(\fq;M)=e(x_1, x_2, \ldots, x_{d_k}; M_k)+\sum_{j=1}^{d_k}\binom{d_k-1}{j-1}\ell(H^j_\fm(M/M_k)),$$
and 
$$e_{d-i}(\fq;M)=\sum_{j=1}^i\binom{i-1}{j-1}\ell(H_\fm^j(M/M_k))$$
for $d_k<i<d_{k+1}$.

Now, let $d_1=1$. We have 
$$\ell(M/\fq^{n+1}M)=\ell(M/\fq^{n+1}M+M_1)+\ell(M_1/\fq^{n+1}M\cap M_1).$$
By Artin-Rees Lemma and the fact that $(x_2, \ldots, x_d)M_1=0$, there is an $n_0>0$ such that $\fq^{n+1}M\cap M_1=\fq^{n+1-n_0}(\fq^{n_0}M\cap M_1)=x_1^{n+1-n_0}(\fq^{n_0}M\cap M_1)$ for all $n+1\geqslant n_0$. Hence,
\begin{multline*}
\ell(M/\fq^{n+1}M)=\ell(M/\fq^{n+1}M+M_1)+\ell(M_1/\fq^{n_0}M\cap M_1)\\+\ell((\fq^{n_0}M\cap M_1)/x_1^{n+1-n_0}(\fq^{n_0}M\cap M_1)).
\end{multline*}
This implies that 
$$e_{d-1}(\fq;M)=e_{d-1}(\fq; M/M_1)+e(x_1;\fq^{n_0}M\cap M_1)=e_{d-1}(\fq; M/M_1)+e(x_1;M_1)$$
 and $e_{d-i}(\fq; M)=e_{d-i}(\fq; M/M_1)$ for all $i>1$. Observe that $M/M_1$ has a generalized Cohen-Macaulay filtration $0\subset M_2/M_1\subset \ldots \subset M_t/M_1=M/M_1$ with $\dim M_2/M_1=d_2 >1$. Then applying the previous argument for $d_1>1$ to the module $M/M_1$ we get the conclusion.
\end{proof}
\begin{corollary}
Keep all notations and hypotheses in Theorem \ref{hs6} . Then the difference
$$\ell(M/\fq^{n+1}M)-\sum_{k=1}^t\binom{n+d_k}{d_k}e(x_1, \ldots, x_{d_k}; M_k)=I_n(M)$$
is independent of the choice of systems of parameters, which are dd-sequences of $M$,  and of the generalized Cohen-Macaulay filtrations of $M$. Moreover, 
$$I_n(M)=\sum_{k=0}^{t-1}\sum_{i=d_k}^{d_{k+1}-1}\binom{n+i}{i}\sum_{j=1}^{i}\binom{i-1}{j-1}\ell(H_\fm^j(M/D_k))+\ell(H^0_\fm(M)).$$
\end{corollary}
\begin{proof}
It is clear from Theorem \ref{hs6} that
$$I_n(M)=\sum_{k=0}^{t-1}\sum_{i=d_k}^{d_{k+1}-1}\binom{n+i}{i}\sum_{j=1}^{i}\binom{i-1}{j-1}\ell(H_\fm^j(M/M_k))+\ell(H^0_\fm(M)).$$
Let $\D: D_0\subset D_1\subset \ldots \subset D_t=M$ be the dimension filtration of $M$. By Lemma \ref{cuongnhan}, $D_i/M_i$ is of finite length for $i=0, 1, \ldots, t$. Hence $H^j_\fm(M/M_i)\simeq H^j_\fm(M/D_i)$ for all $j>0$ and 
$$I_n(M)=\sum_{k=0}^{t-1}\sum_{i=d_k}^{d_{k+1}-1}\binom{n+i}{i}\sum_{j=1}^{i}\binom{i-1}{j-1}\ell(H_\fm^j(M/D_k))+\ell(H^0_\fm(M))$$
does not depend on the system of parameters $\un x$ and the filtration $\F$.
\end{proof}
The following immediate consequence of Theorem \ref{hs6} is a well-known result in the theory of generalized Cohen-Macaulay modules (see \cite{sv}).
\begin{corollary}
Let $M$ be a generalized Cohen-Macaulay module and $\un x=(x_1, \ldots, x_d)$ a standard system of parameters of $M$. Set $\fq=(x_1, \ldots, x_d)$. Then
$$\ell(M/\fq^{n+1}M)=\binom{n+d}{d}e(\un x; M)+\sum_{i=1}^{d-1}\binom{n+i}{i}\sum_{j=1}^{i}\binom{i-1}{j-1}\ell(H_\fm^j(M))+\ell(H^0_\fm(M)).$$
Moreover, the difference
$$\ell(M/\fq^{n+1}M)-\binom{n+d}{d}e(\un x; M)=\sum_{i=1}^{d-1}\binom{n+i}{i}\sum_{j=1}^{i}\binom{i-1}{j-1}\ell(H_\fm^j(M))+\ell(H^0_\fm(M)),$$
is independent of the choice of the standard systems of parameters $\un x$.
\end{corollary}

\end{document}